\documentclass[11pt,a4paper]{article}

\usepackage[utf8]{inputenc}
\usepackage[T1]{fontenc}
\usepackage{lmodern}
\usepackage[margin=2.6cm]{geometry}
\usepackage{amsmath,amssymb,amsthm,mathtools}
\usepackage{booktabs,array}
\usepackage{graphicx}
\graphicspath{{figures/}}
\usepackage{enumitem}
\usepackage[numbers,sort&compress]{natbib}
\usepackage[colorlinks=true,linkcolor=blue!50!black,citecolor=green!40!black,urlcolor=blue!50!black]{hyperref}
\usepackage{xcolor}

\theoremstyle{plain}
\newtheorem{theorem}{Theorem}[section]
\newtheorem{proposition}[theorem]{Proposition}
\newtheorem{corollary}[theorem]{Corollary}
\newtheorem{lemma}[theorem]{Lemma}
\newtheorem{conjecture}[theorem]{Conjecture}
\theoremstyle{definition}
\newtheorem{definition}[theorem]{Definition}
\theoremstyle{remark}
\newtheorem{remark}[theorem]{Remark}

\newcommand{\R}{\mathbb{R}}
\newcommand{\eps}{\varepsilon}
\newcommand{\Gf}{G_0}
\newcommand{\Gs}{G_1}
\newcommand{\ce}[1]{\mathrm{#1}}
\newcommand{\one}{\mathbf 1}
\newcommand{\Net}[1]{\mathcal N_{#1}}

\title{Siphon faces in slow--fast mass-action networks}
\author{Miruna-Mihaela Beldiman\thanks{Institute of Mathematical Statistics
and Applied Mathematics of the Romanian Academy, Bucharest, Romania.
ORCID 0000-0001-9947-7594.}
\and
Andrei-Dan Halanay\thanks{University of Bucharest, Bucharest, Romania.
ORCID 0009-0005-5481-8021. Corresponding author:
\href{mailto:andrei.halanay@unibuc.ro}{andrei.halanay@unibuc.ro}.}}
\date{}

\begin{document}
\maketitle

\begin{abstract}
When the species of a siphon in a mass-action network, however many, are fast and
their invariant face loses transverse stability along the slow flow, trajectories
follow the face past this loss before leaving it. We study three consequences.
(i) If the transverse Metzler linearisation is irreducible, the entry--exit
relation with its Perron root gives the exit. Under acyclic feeding between its
strongly connected components, a max-plus recursion over the reactions gives it;
in general neither the Perron root nor the diagonal blocks do. (ii) With one slow
species, acyclic feeding and a transversal exit, a fold of the interior critical
manifold followed by a drop onto the face gives a unique, exponentially attracting
relaxation oscillation whose turning points tend to the fold and exit point.
Exactly $3474$ minimal planar networks of molecularity $\le3$, with fast
per-capita rate affine and decreasing in the slow species, meet our sufficient
conditions for it. (iii) A siphon-breaking reaction of rate $\delta$ and finite
min-plus breaking order $w$ unfolds the transcritical point into a fold if a
Perron-weighted second-order self-interaction is positive, and not if
negative; the fold's critical component lies at distance of order $\delta^{w/2}$ from the
face. For $\delta=e^{-\mu/\eps}$ the broken face is a floor at $-\mu$ in the
recursion, and planar cycles leave it between the entry--exit and
transcritical points, depending on $\mu$. Simulations confirm predicted exits
and cycles, including cyclic feeding, where the prediction is conjectural, and
high-precision computations confirm the folds; the errors shrink as $\eps$ or
$\delta$ decreases.
\end{abstract}


(secondary).

%======================================================================
\section{Introduction}\label{sec:intro}
%======================================================================

A \emph{siphon} of a reaction network is a set $\Sigma$ of species such that every
reaction producing a species of $\Sigma$ has a species of $\Sigma$ among its
reactants \citep{AngeliDeLeenheerSontag2007,ShiuSturmfels2010}. Under mass-action
kinetics the face $\{x_\Sigma=0\}$ of the positive orthant is then invariant.
Siphons are the standard tool for persistence: in a conservative network in which
every siphon contains the support of a nonnegative conservation law, no species can
be lost \citep{AngeliDeLeenheerSontag2007}. Here we look at the
opposite situation, in a slow--fast network whose fast species include a siphon.
The critical manifold, the set of equilibria of the fast subsystem, then contains the face, and the face is normally
attracting or repelling according to the sign of the leading eigenvalue of the
transverse linearisation. When the slow flow on the face carries the trajectory
across a point where this sign changes, the trajectory does not leave at once. It
follows the face, now repelling, until the accumulated expansion compensates the
accumulated contraction, and only then leaves: for one transverse variable with
eigenvalue $\lambda(\tau)$ along the slow flow, it enters at $\tau_{\rm in}$ and leaves at
the first $\tau_{\rm out}$ with $\int_{\tau_{\rm in}}^{\tau_{\rm out}}\lambda\,d\tau=0$. This is the delayed loss of
stability described by the entry--exit function
\citep{Benoit1981,Schecter1985,DeMaesschalckSchecter2016}; see \citet{Neishtadt1987}
for its counterpart at a Hopf bifurcation. At the point where the face changes stability,
an interior branch of the critical manifold usually crosses the face, in a
transcritical singularity \citep{KrupaSzmolyan2001trans}.

In our earlier papers \citep{BeldimanHalanay2026,BeldimanHalanay2026rankone} we
studied canards, trajectories that follow a repelling branch of the critical
manifold, at interior folds of mass-action networks. The boundary passage is
the other mechanism by which slow--fast networks produce large oscillations. It
appears in predator--prey systems \citep{Hsu2019,HsuRuan2021,AiSadhu2020,LiLiWang2023}
and in enzyme models such as Olsen's \citep{Olsen1983,KuehnSzmolyan2015}.
Olsen's model and the settings of \citet{Hsu2019} and \citet{HsuRuan2021} lie
outside our assumptions (Section~\ref{sec:eenum} and Remark~\ref{rem:landing}). Most of
the entry--exit theory treats one transverse variable, or several that vanish on separate
faces; recent work on coupled transverse variables is discussed in
Section~\ref{sec:disc}. Networks are different in three ways:
\begin{itemize}[nosep]
  \item a siphon usually has several species, and its transverse linearisation is
        a Metzler matrix, often reducible;
  \item its species are coupled not only through this matrix, but also through
        slow reactions and through reactions with two siphon molecules among their
        reactants, which vanish to higher order on the face;
  \item a siphon is a structural property, and it can be destroyed by a single
        reaction with a small rate.
\end{itemize}
Each of the three can move the exit point, which is a turning point of the
oscillations through the face (Sections~\ref{sec:ee}, \ref{sec:ro}
and~\ref{sec:floor}).

\paragraph{Results.} Section~\ref{sec:setting} fixes the setting: the siphon is a
siphon of the whole network, so its face is invariant for $\eps>0$, and no fast
reaction changes other species.

Section~\ref{sec:ee} is about entry--exit. If the transverse Metzler matrix
$\Lambda$ is irreducible, the exit point is given by the classical relation with
the Perron root of $\Lambda$, to $O(\eps)$ (Theorem~\ref{thm:perron}). The proof
is short: the left Perron vector is a linear functional on the positive cone, and
no adiabatic argument is needed. If $\Lambda$ is reducible, the logarithms of the
concentrations of its strongly connected components follow a max-plus recursion
over all reactions that feed one component from another, slow ones and those with
two siphon molecules among their reactants included (Theorem~\ref{thm:tropical},
for an acyclic feeding relation and a transversal exit). The Perron integral gives
the earliest possible exit and the diagonal blocks the latest
(Corollary~\ref{cor:bounds}); both can be wrong.

Section~\ref{sec:ro} uses the exit point to build relaxation oscillations. A
\emph{fold--face cycle} slides along an attracting interior branch to a fold,
drops onto the face, follows the face to the exit point and rises back to the
branch. For one slow species, no fast species outside the siphon, acyclic feeding
and a transversal exit, it perturbs to a unique, exponentially attracting periodic
orbit (Theorem~\ref{thm:ro}). The fold fixes the entry point, so no selection or
Floquet-type condition is needed, unlike for the cycles of \citet{Hsu2019} and
\citet{HsuRuan2021}. The key estimate is that the passage along the face has
derivative $O(1/\eps)$ (Lemma~\ref{lem:face}).

Section~\ref{sec:planar} turns to planar networks, with one fast and one slow
species. There five conditions on the reactions are sufficient for a fold--face
cycle (Proposition~\ref{prop:planar}). Among the networks of molecularity at most
three of \citet{BeldimanHalanay2026} whose fast per-capita rate is affine in the
slow species with negative slope, there are exactly $3474$ minimal networks
satisfying these conditions for some rate constants (Theorem~\ref{thm:enum}),
identified by an exact locking certificate. Each contains a reaction with three
molecules of $\ce X$ among its reactants, so none is bimolecular. Outside this
class there are others (Remark~\ref{rem:outside}).

Section~\ref{sec:break} breaks the siphon. Adding a reaction that produces siphon
species from outside the siphon, at rate $\delta$, unfolds the transcritical
point. If the min-plus \emph{breaking order} $w$ is finite, a fold appears if
$c_2=\ell\cdot Q(r,r)>0$, and not if $c_2<0$; here $c_2$ is a second-order
self-interaction weighted by the left and right null vectors of $\Lambda$. The
critical component of the fold lies at distance of order $\delta^{w/2}$ from the
face (Theorem~\ref{thm:break}). In the planar class $c_2>0$ says exactly that the
interior branch leaves the transcritical point increasing.

Section~\ref{sec:floor} makes the breaking exponentially weak. If
$\delta=e^{-\mu/\eps}$, the breaking reaction acts in the recursion as a feed of
constant logarithmic size $-\mu$, so the broken face is a floor at $-\mu$
(Proposition~\ref{prop:floor}). For planar networks the
lower turning point of the cycle moves from the transcritical point to the
entry--exit point as $\mu$ increases from $0$ to the depth of the passage
(Corollary~\ref{cor:floor}). So for these networks the joint limit
$\delta,\eps\to0$ depends on $\mu=\eps\log(1/\delta)$.

Each of Sections~\ref{sec:ee}--\ref{sec:floor} ends with a numerical check.
Section~\ref{sec:disc} discusses the
literature and open problems. The main open problem is the case of cyclic feeding
(Conjecture~\ref{conj:cyclic}).

%======================================================================
\section{Setting}\label{sec:setting}
%======================================================================

We use slow--fast mass-action networks as in
\citet{BeldimanHalanay2026rankone}, in slow time $\tau$:
\[
  \eps\,x'=N_0v(x)+\eps N_1v(x),
\]
with fast subnetwork $\Gf$ (stoichiometric matrix $N_0$) and slow subnetwork
$\Gs$; $v$ is the vector of mass-action rates. Let $\Sigma$ be a set of $k$
species. Write $z=x_\Sigma\in\R^k_{\ge0}$ and $u$ for the other concentrations. We
assume:
\begin{enumerate}[label=(S\arabic*),nosep]
  \item $\Sigma$ is a siphon of the \emph{whole} network $\Gf\cup\Gs$. So the face
        $H_\Sigma=\{z=0\}$ is invariant for every $\eps\ge0$.
  \item No reaction of $\Gf$ changes a species outside $\Sigma$. So
        $u'=f(u,z,\eps)$ with $f$ smooth. On the face the slow flow is
        $u'=f(u,0,0)$, given by the reactions of $\Gs$ with no species of
        $\Sigma$ among their reactants.
\end{enumerate}
If $\Sigma$ is a siphon of $\Gf$ but not of $\Gs$, or of neither, the face is not
invariant; this is the siphon breaking of Sections~\ref{sec:break}
and~\ref{sec:floor}. Remark~\ref{rem:landing} discusses (S2).

Every term of $(N_0v)_i$ and $(N_1v)_i$, $i\in\Sigma$, has a species of $\Sigma$
among its reactants. Hence
\[
  \eps z'=M(u,z,\eps)\,z,\qquad M(u,0,0)=\Lambda(u).
\]
We call $\Lambda(u)$ the \emph{siphon matrix}. Its entries come from the reactions
of $\Gf$ with exactly one molecule of $\Sigma$ among their reactants, and its
off-diagonal entries are nonnegative: it is a Metzler matrix. Its off-diagonal
sign pattern does not depend on $u>0$, since every off-diagonal entry is a sum of
monomials in $u$ with positive coefficients; so neither does irreducibility. The
diagonal entries are differences, and their change of sign along the slow flow is
what produces the exit. For a Metzler matrix $A$ we write $s(A)$ for its Perron
root, the eigenvalue of largest real part, which is real
\citep[see e.g.][]{Feinberg2019}.

\begin{definition}\label{def:graph}
The \emph{siphon graph} has the species of $\Sigma$ as nodes and an edge
$X_j\to X_i$ ($i\neq j$) if a reaction with exactly one molecule of $\Sigma$, namely
$X_j$, among its reactants has positive net production of $X_i$. The edge is
\emph{fast} or \emph{slow} with the reaction. The \emph{components} are the strongly
connected components of the fast edges, so each diagonal block $\Lambda_{CC}$ is
irreducible. A \emph{feed} of a component $C$ is a reaction $r$ (fast or slow) with
positive net production of a species of $C$ whose reactant vector $y_r$ has
$y_r^\Sigma\ne0$ and no siphon reactant in $C$. Its \emph{weight} is the vector
$y^\Sigma_r$. The \emph{feeding relation} on
components is $D\Rightarrow C$ if $D\ne C$ and some feed of $C$ has a siphon
reactant in $D$.
\end{definition}

The feeds of $C$ are the edges into $C$ from other components, of either speed,
and the reactions with at least two siphon molecules among their reactants, none of
them in $C$, that produce a species of $C$ (``hyperedges''). A reaction with a
siphon reactant $X_j\in C$ is not a feed: it is either linear and fast, and part of
$\Lambda_{CC}$, or linear and slow, of rate at most $\eps Kz_j$, or nonlinear, of
rate at most $K|z|_1z_j$. On the logarithmic scale of Definition~\ref{def:trop} it would offer
at most $\psi_C$, the level of $C$ itself, so it never matters.

Fix an orbit $u_0(\tau)$, $\tau\in[0,T]$, of the slow flow on the face, with
$u_0(0)=u_{\rm in}$. Write $s_C(\tau)=s(\Lambda_{CC}(u_0(\tau)))$ for the Perron
root of block $C$, $s(\tau)=s(\Lambda(u_0(\tau)))=\max_Cs_C(\tau)$, and
$\Psi_C(\sigma,\tau)=\int_\sigma^\tau s_C$. We assume the face is attracting at
entry, $s(0)<0$. Fix a small $\rho>0$. A trajectory \emph{enters} at
$(u_{\rm in},z(0))$ with $z(0)\ge0$, $|z(0)|_1\le\rho$, and its \emph{exit time}
$\tau_\eps$ is the first $\tau\ge\tau_0$ with $|z(\tau)|_1=\rho$, for a fixed small
$\tau_0>0$. Theorem~\ref{thm:perron} takes entries of size exactly $\rho$; smaller
entries are covered by Theorem~\ref{thm:tropical}, through their entry values.

%======================================================================
\section{Entry--exit at the face of a siphon}\label{sec:ee}
%======================================================================

The entry--exit relation rests on a short computation. Take one transverse
variable $z\ge0$ whose face $\{z=0\}$ is invariant:
\[
  \eps z'=z\,g(u,z,\eps),\qquad u'=f(u,z,\eps),\qquad g(u,0,0)=\lambda(u).
\]
A trajectory with $z(0)>0$ keeps $z>0$, and $\phi=\eps\log z$ satisfies
\[
  \phi'=\lambda(u)+O(z+\eps).
\]
While $z$ is small, $u$ stays close to the slow orbit $u_0$ on the face, and $\phi$
follows $\phi(0)+\int_0^\tau\lambda(u_0)$. A trajectory that enters at a distance of
order one, $\phi(0)\approx0$, therefore stays exponentially close to the face as long
as this integral is negative, also after $\lambda$ has become positive, and it leaves
when the integral returns to $0$. The computation uses three facts. The face is
invariant, so the equation can be divided by $z$. The distance from the face is a
single positive number. And on the scale $\eps\log$ an exponentially small quantity
$e^{-c/\eps}$ has the finite size $-c$, while factors of order $\eps$ or $1$ have
size $0$, so the relation is insensitive to prefactors.

In a mass-action network the first fact is assumption (S1), and the positivity of
$z$ needs no assumption. The second fails as soon as the siphon has more than one
species: $z$ is then a vector, whose components are coupled through $\Lambda$,
through slow reactions and through reactions with several siphon molecules among
their reactants. The results of this section start from the question of what
replaces the single number $\eps\log z$. If $\Lambda$ is irreducible, the left
Perron vector $\ell>0$ provides one: the linear functional $\ell\cdot z$ controls
$|z|_1$ because $z\ge0$, and the computation above goes through with the Perron root
$s$ in place of $\lambda$ (Theorem~\ref{thm:perron}). If $\Lambda$ is reducible,
there is one level $\psi_C\approx\eps\log(\ell_C\cdot z_C)$ for each strongly
connected component $C$. Each level grows at the rate of its own block, but its
component is also produced by reactions from other components. On the logarithmic
scale a sum becomes a maximum, $\eps\log(e^{a/\eps}+e^{b/\eps})=\max(a,b)+O(\eps)$,
and a product becomes a sum. So the feeds enter through a maximum, a feed with siphon
reactant vector $y^\Sigma$ contributes $\sum_jy_j\psi_{c(j)}$, and the single integral
of the classical relation becomes the max-plus recursion of
Definition~\ref{def:trop} (Theorem~\ref{thm:tropical}).

The logarithmic scale also organises the rest of the paper. The planar networks of
Section~\ref{sec:planar} have a one-species siphon, and the classical relation fixes
their exit point. The relaxation oscillations of Section~\ref{sec:ro} need the
derivative of the passage along the face, not only its exit point. In the
coordinates $\log z_i$ this derivative is $O(1/\eps)$ (Lemma~\ref{lem:face}), a
loss that the exponential contraction along the attracting interior branch
absorbs. In Section~\ref{sec:break} a small rate $\delta$ of a reaction that breaks
the siphon plays the role of $e^{-1/\eps}$, and the breaking order is the min-plus
counterpart of the recursion. In Section~\ref{sec:floor}, where
$\delta=e^{-\mu/\eps}$, the breaking reaction has logarithmic size $-\mu$ and enters
the recursion as one more feed. The logarithmic scale ignores prefactors, and the
proof of Theorem~\ref{thm:ro} needs bounds on the original scale; the a priori
bounds of Section~\ref{sec:apriori} supply them.

\subsection{Irreducible siphon matrix}

If the siphon matrix is irreducible, the classical entry--exit relation holds, with
the Perron root $s$ of $\Lambda$ in place of the single transverse eigenvalue. A
trajectory that enters near the face while it is attracting ($s<0$) stays near it
after $s$ changes sign. It leaves at the time $\tau^*$ at which contraction and
expansion balance, $\int_0^{\tau^*}s=0$, up to $O(\eps)$.

\begin{theorem}[Perron entry--exit]\label{thm:perron}
Assume (S1), (S2), and that $\Lambda(u)$ is irreducible (for one, hence every,
$u>0$). Let $\Psi(\tau)=\int_0^\tau s$, and suppose there is
$\tau^*\in(0,T)$ with $\Psi<0$ on $(0,\tau^*)$, $\Psi(\tau^*)=0$ and $s(\tau^*)>0$.
Then there are $\rho_0,C>0$ such that for $\rho\le\rho_0$, small $\eps>0$ and
every entry $z(0)\ge0$ with $|z(0)|_1=\rho$,
\[
  |\tau_\eps-\tau^*|\le C\eps,\qquad |u(\tau_\eps)-u_0(\tau^*)|\le C\eps .
\]
More precisely, with the left Perron vector $\ell(u)>0$ of $\Lambda(u)$,
$\ell\cdot\one=1$,
\begin{equation}\label{eq:phi}
  \eps\log\bigl(\ell(u(\tau))\cdot z(\tau)\bigr)=\eps\log\bigl(\ell(u_{\rm in})\cdot z(0)\bigr)+\Psi(\tau)+O(\eps)
  \qquad\text{uniformly on }[0,\tau_\eps].
\end{equation}
\end{theorem}

\begin{proof}
On a compact neighbourhood $\mathcal K$ of the orbit, $\Lambda(u)$ is irreducible,
so its Perron root $s(\Lambda(u))$ is a simple eigenvalue. Hence $s$ and the normalised left
eigenvector $\ell(u)$ are smooth, and $\ell(u)\ge c\one$ for some $c>0$. For $z\ge0$,
$|z|_1\le(\ell\cdot z)/c$.

Let $Z=\ell(u)\cdot z$. Using $\ell\Lambda=s\ell$,
\[
  \eps Z'=\ell M z+\eps(D\ell\,u')\cdot z
        =s(u)Z+\ell\bigl(M(u,z,\eps)-\Lambda(u)\bigr)z+\eps(D\ell f)\cdot z .
\]
The last two terms are $O(|z|_1+\eps)\,|z|_1=O(|z|_1+\eps)\,Z$, by positivity
of $z$. So $\phi=\eps\log Z$ satisfies
\begin{equation}\label{eq:phiprime}
  \phi'=s(u(\tau))+O(|z(\tau)|_1+\eps).
\end{equation}
Since
$u'=f(u,0,0)+O(|z|_1+\eps)$, Gronwall gives
$|u(\tau)-u_0(\tau)|\le C(\eps+\int_0^\tau|z|_1)$, and so
$s(u(\tau))=s(\tau)+O(\eps+\int_0^\tau|z|_1)$.

It remains to show $\int_0^{\tau_\eps}|z|_1=O(\rho\eps)$. Choose $\tau_1>0$ with
$s<s(0)/2$ on $[0,\tau_1]$, and $\eta>0$ with $s>s(\tau^*)/2$ on
$[\tau^*-\eta,\tau^*+\eta]$ and $\tau^*+\eta<T$. Let
$m=-\max_{[\tau_1,\tau^*-\eta]}\Psi>0$. Fix $M>0$, chosen below. As long as
$\int_0^\tau|z|_1\le M\rho\eps$ and $\rho$, $\eps$ are small, \eqref{eq:phiprime}
gives:
\begin{itemize}[nosep]
  \item on $[0,\tau_1]$, $\phi'\le s(0)/4$, so $|z|_1\le C\rho e^{s(0)\tau/(4\eps)}$
        and the integral over $[0,\tau_1]$ is at most $C\rho\eps$;
  \item on $[\tau_1,\tau^*-\eta]$, $\phi\le\eps\log\rho-m/2$, so
        $|z|_1\le Ce^{-m/(2\eps)}$;
  \item on $[\tau^*-\eta,\tau_\eps]$, $\phi'\ge s(\tau^*)/4$ and $|z|_1\le\rho$, so
        the integral is at most $C\rho\eps$ again, by integrating backwards from
        $\tau_\eps$.
\end{itemize}
Run these estimates on $[0,\min(\tau_\eps,\tau^*+\eta)]$. They give
$\int|z|_1\le C_0\rho\eps$ there, with $C_0$ independent of $M$ once $\eps$ is small
(depending on $M$). Taking $M=2C_0$, a continuity argument closes the bootstrap on
this interval. Then $\tau_\eps\le\tau^*+\eta$: otherwise, by \eqref{eq:phiprime},
$\phi(\tau^*+\eta)\ge\eps\log(c\rho)+\Psi(\tau^*+\eta)-C\eps>\eps\log\rho$ for small
$\eps$, while $\ell\cdot z\le|z|_1\le\rho$ before $\tau_\eps$. Integrating \eqref{eq:phiprime} now
gives \eqref{eq:phi}. At $\tau_\eps$, $|z|_1=\rho$, so
$\phi(\tau_\eps)-\phi(0)\in[\eps\log c,-\eps\log c]$. Hence
$\Psi(\tau_\eps)=O(\eps)$, and $\tau_\eps=\tau^*+O(\eps)$ because $\Psi$ has no zero
in $(0,\tau^*)$, by the bound on $(0,\tau^*-\eta]$, and crosses zero transversally at
$\tau^*$. The bound on $u$ follows from Gronwall.
\end{proof}

\begin{remark}
\begin{itemize}[nosep]
  \item The proof needs no adiabatic argument (see the start of this section). This
        simplifies the proof, but does not weaken the hypotheses:
        for an irreducible Metzler family over a compact set the Perron root is in
        any case separated from the rest of the spectrum.
  \item Nonlinear terms of $M$ with two or more siphon molecules are harmless
        here, because $\ell>0$ controls all components.
  \item For $k=1$ this is the classical entry--exit relation
        \citep{Benoit1981,Schecter1985,DeMaesschalckSchecter2016}. For $k=2$, when the second eigenvalue stays negative, the
        statement is proved with nonstandard analysis by
        \citet{BacaerLobrySari2020}, for the extinction probability of a two-type
        branching process in a slow periodic environment, where the generating
        function has an invariant face with an irreducible Metzler linearisation.
        They report that they could not prove it for more than two types, and
        give numerics for four.
\end{itemize}
\end{remark}

\subsection{Reducible siphon matrix: a tropical recursion}

\begin{definition}\label{def:trop}
Let the entry values $a_C=\lim_{\eps\to0}\eps\log(\ell_C\cdot z_C(0))\in[-\infty,0]$
be given, with $\ell_C$ the left Perron vector of $\Lambda_{CC}(u_{\rm in})$. The
\emph{tropical recursion} is
\begin{equation}\label{eq:trop}
  \psi_C(\tau)=\max\Bigl\{a_C+\Psi_C(0,\tau),\ \sup_{0\le\sigma\le\tau}\bigl[F_C(\psi(\sigma))+\Psi_C(\sigma,\tau)\bigr]\Bigr\},
\end{equation}
\[
  F_C(\psi)=\max_{\text{feeds }r\text{ of }C}\ \sum_{j\in\Sigma}y_{rj}\,\psi_{c(j)},
\]
where $c(j)$ is the component of $X_j$ (and $F_C=-\infty$ if $C$ has no feed), with
$-\infty+x=-\infty$. The \emph{tropical exit time} is
$\tau_{\rm ex}=\inf\{\tau>0:\max_C\psi_C(\tau)=0\}$.
\end{definition}

Since the siphon reactants of a feed of $C$ lie outside $C$, $F_C$ depends only on
the other components. An entry value $a_C=-\infty$, for instance $z_C(0)=0$, is
allowed; $\psi_C$ is then finite from $\tau=0$ on if some feed offers a finite
value at $\tau=0$, and $-\infty$ throughout otherwise.

In words: $\psi_C$ grows at the rate $s_C$ of its own block, unless a feed
offers more; a feed through a reaction with siphon reactant vector $y^\Sigma$
offers $\sum_jy_j\psi_{c(j)}$, the logarithmic size of its rate. A slow feed
carries a factor $\eps$ in its rate, but $\eps\log\eps\to0$, so on the logarithmic
scale it counts as much as a fast one. Equivalently, as a Skorokhod problem, reflected at the feed level $F_C(\psi)$:
$\psi_C(0)=\max(a_C,F_C(\psi(0)))$, $\psi_C\ge F_C(\psi)$, and
$\psi_C-\Psi_C(0,\cdot)$ is nondecreasing and constant off the contact set
$\{\psi_C=F_C(\psi)\}$. (Without the monotonicity, a solution could stay in contact
with a falling feed after its own block has started to grow.)

\begin{definition}[Transversal exit]\label{def:transversal}
Let $E=\{C:\psi_C(\tau_{\rm ex})=0\}$ be the \emph{exiting} components, and
$E^*\subset E$ the \emph{self-exiting} ones: those $C\in E$ no feed of which has a
siphon reactant in a component of $E$. The exit is \emph{transversal} if
$\tau_{\rm ex}<T$, $s_C(\tau_{\rm ex})>0$ for $C\in E^*$, and $s_C(\tau_{\rm ex})<0$
for $C\in E\setminus E^*$.
\end{definition}

For acyclic feeding $E^*\ne\emptyset$: a component of $E$ that comes first in the
feeding order is self-exiting. Its feeds offer a negative value at $\tau_{\rm ex}$,
so near $\tau_{\rm ex}$ it follows its own block and crosses $0$ with slope $s_C$.
The other components of $E$ are carried to $0$ by their feeders, like $\ce X_2$ in
network R2 below, a species that decays on its own at the exit. Components of $E\setminus E^*$
that grow on their own at the exit, and zero rates at the exit, are excluded for
simplicity. Theorem~\ref{thm:tropical} uses only the condition on $E^*$; the
condition on $E\setminus E^*$ is used for the face passage
(Lemma~\ref{lem:integral} and Section~\ref{sec:ro}).

\begin{theorem}[Tropical entry--exit, acyclic case]\label{thm:tropical}
Assume (S1), (S2), $s(0)<0$, and that the feeding relation is acyclic. Let the entry
values satisfy $\eps\log(\ell_C\cdot z_C(0))\to a_C$, uniformly in the sense that the
error is $o(1)$ as $\eps\to0$ (for $a_C=-\infty$: the left side tends to $-\infty$,
with $\log0=-\infty$), with $\max_Ca_C=0$. Then
\[
  \eps\log\bigl(\ell_C(u(\tau))\cdot z_C(\tau)\bigr)\to\psi_C(\tau)
\]
uniformly on $[\theta,\tau_{\rm ex}-\eta]$ for all $\eta>0$ and $0<\theta<\tau_{\rm ex}-\eta$
(where $\psi_C=-\infty$, in the sense that the left side is eventually below any
constant), and on $[0,\tau_{\rm ex}-\eta]$ if $a_C\ge F_C(\psi(0))$ for every $C$,
for instance if all $a_C=0$. If moreover $\tau_{\rm ex}<T$ and
$s_C(\tau_{\rm ex})>0$ for every $C\in E^*$, in particular if the exit is
transversal, then $\tau_\eps\to\tau_{\rm ex}$ and $u(\tau_\eps)\to u_0(\tau_{\rm ex})$.
\end{theorem}

\begin{proof}
Order the components so that every feed of $C$ has its siphon reactants in earlier
components; this is possible because the feeding relation is acyclic.
Write $Z_C=\ell_C(u)\cdot z_C$ and $\phi_C=\eps\log Z_C$.

\emph{Equation for $Z_C$.} As in the proof of Theorem~\ref{thm:perron},
\[
  \eps Z_C'=\bigl(s_C(u)+O(|z|_1+\eps)\bigr)Z_C+R_C(\tau),\qquad
  R_C=\sum_{\text{feeds }r}\kappa_r(u,\eps)\,z^{y_r^\Sigma}\,\bigl(1+O(|z|_1)\bigr),
\]
where $\kappa_r>0$ is the rate constant times the monomial in $u$, times the
$\ell_C$-weighted net production $\ell_C\cdot n_{r,C}>0$ (a feed consumes nothing in
$C$ and produces some species of $C$, and $\ell_C>0$), and carries a factor $\eps$ for slow feeds. The remaining terms
of $\ell_C(M-\Lambda)z$ have a siphon reactant in $C$, so they are
$O(|z|_1)Z_C$, and are absorbed in the first term.

\emph{Components are internally comparable.} Write the equations of the species of
$C$ as
\begin{equation}\label{eq:metzlersplit}
  \eps z_C'=A_C(\tau)z_C+q_C(\tau).
\end{equation}
Here $q_C\ge0$ collects the production terms whose siphon reactants all lie
outside $C$. Every other term goes into $A_C$:
\begin{itemize}[nosep]
  \item a consumption term of $X_i$ has $y_i\ge1$ and is written as
        $-(\kappa z^{y-e_i})\,z_i$;
  \item a production term with a reactant $X_j\in C$ is written as
        $(\kappa z^{y-e_j})\,z_j$.
\end{itemize}
So $A_C$ is Metzler. On $\{|z|_1\le K\rho\}$ its entries are bounded by some
$K_1$, and $(A_C)_{ij}\ge a_0>0$ on the fast edges $X_j\to X_i$ inside $C$. By the
comparison principle for linear cooperative systems, the propagator of
$\eps\zeta'=A_C\zeta$ over a fast time $1$ satisfies
\[
  \Phi_C(\tau,\tau-\eps)\ \ge\ e^{-K_1}e^{a_0B_C}\ \ge\ c_1\one\one^{\top},
\]
where $B_C$ is the adjacency matrix of the fast edges of $C$, which is
irreducible. Since $q_C\ge0$, $z_C(\tau)\ge\Phi_C(\tau,\tau-\eps)z_C(\tau-\eps)$.
Together with $\ell_C\ge c\one$ and $\ell_C\cdot\one=1$ this gives
\[
  c_1Z_C(\tau-\eps)\le z_j(\tau)\le Z_C(\tau)/c\qquad(j\in C,\ \tau\ge\eps).
\]
So once $\phi_C\to\psi_C$ uniformly, with $\psi_C$ continuous,
$\eps\log z_j=\phi_C+o(1)$ uniformly for $\tau\ge\eps$. On the logarithmic scale,
$z^{y^\Sigma}$ then equals $\sum_jy_j\phi_{c(j)}$ up to $o(1)$. For
$\sigma\in[0,\eps]$ only the upper bound $z_j\le Z_C/c$ is used; it holds for all
$\tau$.

\emph{Variation of constants and Laplace's principle.} The equation for $Z_C$ is
linear in $Z_C$, with the source $R_C\ge0$:
\[
  Z_C(\tau)=e^{(\Psi_C(0,\tau)+o(1))/\eps}Z_C(0)
  +\frac1\eps\int_0^\tau e^{(\Psi_C(\sigma,\tau)+o(1))/\eps}R_C(\sigma)\,d\sigma,
\]
where the $o(1)$ terms use the bound $\int|z|_1=o(1)$, proved by the same
bootstrap as before: all $\phi_C$ are below $\eps\log\rho$, and they are bounded
away from $0$ on compact subsets of $(0,\tau_{\rm ex})$ (by the induction hypothesis
and $s(0)<0$). By induction, the earlier components satisfy
$\phi_D\to\psi_D$ uniformly, so $\eps\log R_C(\sigma)\to F_C(\psi(\sigma))$ uniformly,
since $\eps\log\kappa_r\to0$ for fast and slow feeds alike. Two facts finish the
step:
\begin{itemize}[nosep]
  \item $\eps\log(A+B)=\max(\eps\log A,\eps\log B)+O(\eps)$ for $A,B>0$;
  \item Laplace's principle: if $g_\eps\to g$ uniformly on $[0,\tau]$ and $g$ is
        continuous, then
        \[
          \eps\log\Bigl(\frac1\eps\int_0^\tau e^{g_\eps(\sigma)/\eps}\,d\sigma\Bigr)\to\max_{[0,\tau]} g .
        \]
\end{itemize}
They give $\phi_C\le\psi_C+o(1)$ uniformly on $[0,\tau_{\rm ex}-\eta]$, which is the
upper half of \eqref{eq:trop}. For the lower half, the feeders converge only on
$[\theta',\tau_{\rm ex}-\eta]$; restricting the Laplace integral to $[\theta',\tau]$
gives $\phi_C(\tau)\ge\psi_C(\tau)-\omega(\theta')-o(1)$ for $\tau\ge\theta'$, with
$\omega$ a modulus of continuity of $\psi$, and $\theta'$ is arbitrary. This gives
uniform convergence on $[\theta,\tau_{\rm ex}-\eta]$. If $a_C\ge F_C(\psi(0))$ for
every $C$, then $\psi_C(0)=a_C=\lim\phi_C(0)$, and the first term alone gives the
lower bound near $\tau=0$.

\emph{Exit.} Assume the exit is transversal, and fix a small $\eta>0$. All
$\psi_C<0$ on compact subsets of $(0,\tau_{\rm ex})$, so there is no exit before
$\tau_{\rm ex}-\eta$ once $\eps$ is small. Let $C\in E^*$. Its feeds have their
siphon reactants in components outside $E$. Let $K$ bound $|s_C|$, and
$J=[\tau_{\rm ex}-\eta,\tau_{\rm ex}+(3K/s_0)\eta]$, with $\eta$ so small that
$J\subset[0,T)$, $F_C(\psi)\le-2m<0$ on $J$ for some $m>0$, and $s_C\ge s_0>0$ on $J$.
Then $\psi_C$ follows its own block,
$\psi_C(\tau_{\rm ex}-\eta)=-\Psi_C(\tau_{\rm ex}-\eta,\tau_{\rm ex})\ge-K\eta$. As
long as $\tau\in J$ and $\tau\le\tau_\eps$, so that $|z|_1\le\rho$, the equation for
$Z_C$ and $R_C\ge0$ give $\eps Z_C'\ge(s_C-K\rho-K\eps)Z_C\ge\tfrac{s_0}2Z_C$, for
small $\rho$. Hence
\[
  \eps\log Z_C(\tau)\ \ge\ -K\eta-o(1)+\tfrac{s_0}2\bigl(\tau-\tau_{\rm ex}+\eta\bigr),
\]
which exceeds $\eps\log\rho$ before $\tau_{\rm ex}-\eta+3K\eta/s_0$. So
$\tau_{\rm ex}-\eta\le\tau_\eps\le\tau_{\rm ex}+(3K/s_0)\eta$ for small $\eps$, and
$\eta$ is arbitrary. The bound on $u(\tau_\eps)$ follows by Gronwall: the integral
$\int_0^{\tau_\eps}|z|_1$ is at most $o(1)$ before $\tau_{\rm ex}-\eta$ and at most
$K\rho\eta$ after it, so it tends to $0$, since $\eta$ is arbitrary.
\end{proof}

From the Laplace integrals we expect the error of Theorem~\ref{thm:tropical} to be
of order $\eps\log(1/\eps)$, against $O(\eps)$ in Theorem~\ref{thm:perron}; the two
values of $\eps$ for R1 in Table~\ref{tab:num} are consistent with this, and in R2 the error is about $10^{-5}$ at $\eps=10^{-3}$.

\begin{corollary}[Bounds]\label{cor:bounds}
Under the hypotheses of Theorem~\ref{thm:tropical}, let $\tau_P$ be the first zero of
the Perron integral $\int_0^\tau s$, and $\tau_C$ that of $\int_0^\tau s_C$. Then
\[
  \tau_P\le\tau_{\rm ex}\le\min\{\tau_C:\ a_C=0\}.
\]
\end{corollary}

\begin{proof}
The upper bound holds because $\psi_C\ge a_C+\Psi_C(0,\cdot)$. For the lower bound,
show $\psi_C(\tau)\le\int_0^\tau s$ for all $C$, by induction along the order. Since
all $\psi\le0$ before the exit and every weight has $|y^\Sigma|\ge1$, a feed offers
$F_C(\psi(\sigma))\le\max_D\psi_D(\sigma)\le\int_0^\sigma s$, and
$\Psi_C(\sigma,\tau)\le\int_\sigma^\tau s$.
\end{proof}

So the Perron integral gives the earliest possible exit, and the diagonal blocks the
latest. Both bounds can be strict. In R1 below, the block with the largest rate
feeds the exiting component through a linear edge, and the Perron integral is
exact. In R2 the dominant block cannot feed the component that crosses, and the
Perron integral underestimates the exit. A feed through a hyperedge passes on only
a multiple of the level, $2\psi$ instead of $\psi$, and the Perron integral again
underestimates the exit, although the dominant block is upstream (R6).

\begin{conjecture}[Cyclic feeding]\label{conj:cyclic}
Theorem~\ref{thm:tropical} holds without the acyclicity assumption, with $\psi$ the
least solution of \eqref{eq:trop}.
\end{conjecture}

The right side of \eqref{eq:trop} is monotone in $\psi$, so a least solution exists,
by the Knaster--Tarski theorem.

Cycles in the feeding relation arise in two ways. Both pass through feeds that are
small near the face, and are therefore not in $\Lambda$:
\begin{itemize}[nosep]
  \item a slow edge back to an upstream component, of size $\eps$;
  \item a hyperedge, of size $O(|z|)$.
\end{itemize}
Both create spurious solutions, other than the least one. At $\psi=0$ a loop through a hyperedge is neutral
($2\cdot0=0$), and solutions can stay at $0$; for $\psi<0$ it contracts
($2\psi<\psi$), and we expect the least solution to be the limit, as
$\theta\downarrow0$, of the solutions for entry values $a_C-\theta$. A cycle of
linear slow edges is neutral at every level: every $\psi$ shared by the components
on the cycle with $\psi(0)\ge\max a_C$ and
$\psi(\tau)-\psi(\sigma)\ge\max_C\Psi_C(\sigma,\tau)$ solves \eqref{eq:trop} on the
cycle, for instance $\psi\equiv0$ while the rates are negative. The least solution
grows at the largest of the rates on the cycle. Networks R3 and R4 below
support the conjecture in both cases.

\begin{remark}[Related recursions]\label{rem:tropicalrelated}
Recursions of the same max-plus type describe the sizes of subpopulations in
stochastic models of adaptive dynamics in which a mutation occurs with probability
$K^{-1/\alpha}$, with $K$ the carrying capacity. There the exponent
$\beta_w=\log N_w/\log K$ of a trait $w$, on the time scale $\log K$, converges, up to a stopping time, to a
piecewise affine limit: the maximum of the growth of $w$ at its own fitness and of
the levels $\beta_u-d(u,w)/\alpha$ supplied by mutants from other traits $u$, with
$d$ the distance in the mutation graph, truncated at $0$
\citep[Thm.~2.2]{CoquilleKrautSmadi2021}; see also
\citet{BovierCoquilleSmadi2019,ChampagnatMeleardTran2021} and, for the
deterministic large-population limit, \citet{KrautBovier2019}.
\citet[Thm.~2.4]{EsserKraut2025} allow piecewise constant parameters that change periodically on a
time scale faster than $\log K$, and the limit then involves time-averaged fitnesses. With
$\eps=1/\log K$, a mutation from $u$ to $w$ is a feed of $w$ by $u$ whose rate
constant carries the factor $e^{-d(u,w)/(\alpha\eps)}$. Such a factor lowers the
value offered by the feed by $d(u,w)/\alpha$, as the breaking reaction of
Section~\ref{sec:floor} does (Proposition~\ref{prop:floor}). Our setting differs in
four ways. The components are strongly connected blocks, measured by their left
Perron vectors, rather than single types. The rates $s_C$ vary along the slow flow,
so the growth terms are integrals. Feeds through hyperedges offer multiples
$\sum_jy_j\psi_{c(j)}$. And the system is a deterministic slow--fast system near an
invariant face. For a continuous trait, the logarithm $\eps\log n_\eps$ of a
population density was used by \citet{DiekmannJabinMischlerPerthame2005}; its limit
is described by a constrained Hamilton--Jacobi equation
\citep{BarlesMirrahimiPerthame2009}.
\end{remark}

\subsection{A priori bounds for entries inside the orthant}\label{sec:apriori}

Theorems~\ref{thm:perron} and~\ref{thm:tropical} work on the logarithmic scale. The
proof of Theorem~\ref{thm:ro} also needs bounds on the original scale, for entries
that are not exponentially small.

\paragraph{Standing assumptions.}
\begin{itemize}[nosep]
  \item (S1), (S2) hold, and the feeding relation is acyclic. This includes the
        irreducible case, where there is a single component.
  \item The entry point $u_{\rm in}$ lies in a compact set $U_0$ with
        $s(\Lambda(u_{\rm in}))<0$.
  \item The entry $z(0)$ lies in a compact set $K_0\subset(0,\infty)^k$ with
        $|z(0)|_1\le\rho$.
\end{itemize}
Constants may depend on $\rho$, $U_0$ and $K_0$, but not on $\eps$ or on the
entry. The exceptions are the constants $K$ in bounds of the form $|z|_1\le K\rho$
and $K'$ in Lemma~\ref{lem:integral}(c): they depend only on $U_0$, on lower bounds
for $\varpi$ (defined below), the $\ell_C$ and $s_0$, and on bounds for the rate
functions, not on $\rho$, once $\rho$ is small.

Write $\xi_i=\log z_i$. Then $\eps\xi_i'=(\eps z')_i/z_i$ is a sum of
\emph{terms} $n_{ri}\kappa_r z^{y_r}/z_i$, one for each reaction $r$ with
$n_{ri}\ne0$, where $n_{ri}$ is the net change of $X_i$. A term is a
\emph{production} term if $n_{ri}>0$ and a \emph{consumption} term otherwise.

On $[0,\tau_\eps]$, $|z|_1\le K\rho$. To see this, take
$\varpi=\one^{\top}(-\Lambda(u_{\rm in}))^{-1}$. It is positive, since the inverse of a
Hurwitz Metzler matrix is nonnegative with positive diagonal, and
$\varpi\Lambda(u_{\rm in})=-\one^{\top}$. So $\varpi\cdot z$ decreases while $u$ stays
near $u_{\rm in}$ and $|z|_1$ is small. Afterwards $|z|_1\le\rho$
by the definition of $\tau_\eps$.

\begin{lemma}[A priori bounds]\label{lem:apriori}
There is $K$ such that on every interval $[0,\tau']$ on which $|z|_1\le K_2\rho$, for
a fixed $K_2$ (in particular on $[0,\tau_\eps]$, and on $[0,\tau^*_\eps]$ of
Lemma~\ref{lem:integral}(c)):
\begin{enumerate}[label=(\alph*),nosep]
  \item $z_j\le Kz_i$ whenever $X_i$ and $X_j$ lie in the same component;
  \item every term of every $\eps\xi_i'$ is bounded by $K$ in absolute value. In
        particular $|\eps\xi_i'|\le K$.
\end{enumerate}
\end{lemma}

\begin{proof}
By induction along an order of the components in which every feed of $C$ has its
siphon reactants in earlier components. Assume (b) for all earlier
components, and let $X_i\in C$. The terms of $\eps\xi_i'$ are of four kinds.
\begin{enumerate}[label=(\roman*),nosep]
  \item Consumption terms, $-|n|\kappa z^{y-e_i}$. They are bounded, since
        $|z|_1\le K\rho$.
  \item The production term with $y=e_i$, if any. It is constant in $z$.
  \item Production terms with $y\ne e_i$ and a reactant $X_j\in C$. They equal
        $n\kappa z^{y-e_j}\,(z_j/z_i)$, which is bounded once (a) holds.
  \item Feeds, $\Theta=n\kappa z^y/z_i$ with $y$ supported on earlier components.
\end{enumerate}

\emph{Proof of (a) for $C$.} Use the splitting \eqref{eq:metzlersplit}. Every term of
$q_C$ has a logarithmic derivative
$\eps(\log\kappa z^y)'=\sum_jy_j\eps\xi_j'+O(\eps)$ bounded by some $K_q$, by
the induction hypothesis. So
\[
  q_C(\sigma)\le e^{2K_q}q_C(\sigma')\qquad\text{for }|\sigma-\sigma'|\le2\eps .
\]
With $c_1$ and $K_1$ as in the proof of Theorem~\ref{thm:tropical}, and $\Phi_C$
also bounded above by $C_1\one\one^{\top}$ over fast times at most $1$, we get for
$\tau\ge2\eps$ and $X_j\in C$:
\[
  z_j(\tau)\le C_1\one^{\top}z_C(\tau-\eps)
  +\frac{C_1}{\eps}\int_{\tau-\eps}^{\tau}\one^{\top}q_C .
\]
For every $X_i\in C$,
\[
  z_i(\tau)\ge c_1\one^{\top}z_C(\tau-\eps),\qquad
  z_C(\tau-\eps)\ge\frac{e^{-K_1}}{\eps}\int_{\tau-2\eps}^{\tau-\eps}q_C .
\]
Hence
\[
  z_i(\tau)\ge\tfrac{c_1}2\one^{\top}z_C(\tau-\eps)+\tfrac{c_1}{2}e^{-K_1-2K_q}\,
  \frac1\eps\int_{\tau-\eps}^{\tau}\one^{\top}q_C
  \ \ge\ K^{-1}z_j(\tau).
\]
For $\tau\le2\eps$:
\begin{itemize}[nosep]
  \item $z_i(\tau)\ge e^{-2K_1}z_i(0)\ge e^{-2K_1}\min_{K_0}z_i$;
  \item $z_j(\tau)\le K\rho$.
\end{itemize}

\emph{Proof of (b) for $C$.} Terms (i)--(iii) are now bounded. For a feed $\Theta$,
every other term of $\eps\xi_i'$ is either a production term ($\ge0$) or a
consumption term ($\ge-K$). So $\eps\xi_i'\ge\Theta-K'$, and
\[
  \eps(\log\Theta)'=\sum_jy_j\eps\xi_j'-\eps\xi_i'+O(\eps)\le K''-\Theta .
\]
Hence $\Theta\le\max\bigl(\Theta(0),K''\bigr)$, and $\Theta(0)$ is bounded on $K_0$.
\end{proof}

\begin{lemma}[Integral bound]\label{lem:integral}
Assume moreover that $U_0$ is a small compact neighbourhood of a point $u_F$ from
which the exit is transversal (Definition~\ref{def:transversal}), and let $E$ and
$E^*$ be the sets of that exit. By continuity, for every $u_{\rm in}\in U_0$:
$\tau_{\rm ex}<T$; near $\tau_{\rm ex}(u_{\rm in})$, $s_C\ge s_0>0$ for $C\in E^*$
and $s_D\le-s_0$ for $D\in E\setminus E^*$; and $\psi_D(\tau_{\rm ex})\le-m'<0$ for
$D\notin E$. Only these facts are used below. All $a_C=0$ here, since $\eps\log z_i(0)\to0$ uniformly on $K_0$. Let
$G_*=\sum_{C\in E^*}\ell_C(u)\cdot z_C$, and let $\tau^*_\eps$ be the first
$\tau\ge\tau_0$ with $G_*=\rho$. Then
\begin{enumerate}[label=(\alph*),nosep]
  \item $\tau_\eps\to\tau_{\rm ex}(u_{\rm in})$ uniformly on $U_0\times K_0$;
  \item $\int_0^{\tau_\eps}|z(\tau)|_1\,d\tau\le K\eps$;
  \item $\tau_\eps\le\tau^*_\eps\le\tau_\eps+K\eps$, $|z|_1\le K\rho$ on
        $[0,\tau^*_\eps]$, and~(b) holds with $\tau^*_\eps$ in place of $\tau_\eps$.
\end{enumerate}
\end{lemma}

\begin{proof}
\emph{Uniformity.} The proof of Theorem~\ref{thm:tropical} uses the entry only
through $\eps\log(\ell_C\cdot z_C(0))$, which tends to $0$ uniformly on $K_0$, and
through $u_{\rm in}$. Every $o(1)$ in that proof depends only on this rate, on the
moduli of continuity of the functions $\Psi_C$ and on
$\sup|s_C(u)-s_C(u_0(\tau))|$ over a compact neighbourhood of the orbits, and these
are uniform for $u_{\rm in}\in U_0$. The functions $\psi$ are equicontinuous in
$\tau$, and depend continuously on $u_{\rm in}$ in the supremum norm: they are
built by finitely many maxima, suprema and sums from the $\Psi_C$, which depend
continuously on $u_{\rm in}$. The function $\max_C\psi_C$ is negative on compact
subsets of $(0,\tau_{\rm ex})$, and positive just after $\tau_{\rm ex}$, where the
components of $E^*$ grow; so $\tau_{\rm ex}$ is continuous in $u_{\rm in}$. A compactness
argument then gives~(a).

\emph{Integral.} Split $[0,\tau_\eps]$ into three windows.
\begin{itemize}[nosep]
  \item \emph{Entry window $[0,\tau_1]$.} With $\varpi$ as above,
        $\varpi\Lambda(u)\le-c\,\varpi$ for $u$ near $U_0$. So
        $\eps(\varpi\cdot z)'\le(-c+K|z|_1+K\eps)\,\varpi\cdot z\le-\tfrac c2\,\varpi\cdot z$.
        Hence $|z|_1\le K\rho e^{-c\tau/(2\eps)}$, and the integral over this window
        is $O(\eps)$.
  \item \emph{Middle window $[\tau_1,\tau_{\rm ex}-\eta]$.} By
        Theorem~\ref{thm:tropical} and (a), $\max_C\phi_C\le-m/2$ there, for some
        $m>0$. So $|z|_1\le Ke^{-m/(2\eps)}$.
  \item \emph{Exit window $[\tau_{\rm ex}-\eta,\tau_\eps]$.} Take $\eta$ small, so
        that on the window $s_C\ge s_0>0$ for $C\in E^*$ and $s_D\le-s_0$ for
        $D\in E\setminus E^*$.
        \begin{itemize}[nosep]
          \item For $D\notin E$, $\psi_D\le-m'<0$ on the window. The upper half of
                the variation-of-constants argument needs no continuity. It gives
                $\phi_D\le\psi_D+K\eta\rho+o(1)$ up to $\tau_\eps$, which is at
                most $-m'/2$. So $|z_D|_1\le Ke^{-m'/(2\eps)}$.
          \item For $C\in E^*$, feeds are nonnegative, so
                $\eps Z_C'\ge(s_C-K|z|_1-K\eps)Z_C\ge\tfrac{s_0}2Z_C$, and
                $Z_C(\tau_\eps)\le\rho$. So
                $Z_C(\tau)\le\rho e^{-s_0(\tau_\eps-\tau)/(2\eps)}$, and its
                integral is at most $2\rho\eps/s_0$.
          \item For $D\in E\setminus E^*$, along the order,
                $\eps Z_D'\le-\tfrac{s_0}2Z_D+R_D$ with
                $R_D\le K\sum_{C}Z_C+Ke^{-m'/(2\eps)}$, the sum over the
                components $C\in E$ that feed $D$: a linear feed is at most $KZ_C$,
                a hyperedge at most $K|z|_1^{|y|-1}Z_C$, and a feed with a siphon
                reactant outside $E$ carries its exponentially small factor.
                Integrating over the window,
                $\tfrac{s_0}2\int Z_D\le\eps Z_D(\tau_{\rm ex}-\eta)+\int R_D$,
                which is $O(\eps)$ by induction.
        \end{itemize}
        All contributions to the integral are $O(\eps)$.
\end{itemize}

\emph{Exit section.} This is a bootstrap on the part of the window where
$|z|_1\le2K'\rho$, with $K'$ fixed below; $\rho$ is small enough that all the
inequalities above hold there. On it, $G_*$ is nondecreasing, since
$\eps G_*'\ge\tfrac{s_0}2G_*$. For $D\in E\setminus E^*$ we have
$Z_D(\tau_{\rm ex}-\eta)\le e^{-c/\eps}$, since $\psi_D(\tau_{\rm ex}-\eta)<0$, and
comparison in $\eps Z_D'\le-\tfrac{s_0}2Z_D+R_D$ gives, by induction along the
order,
\[
  Z_D(\tau)\le K\max_{\sigma\le\tau}G_*(\sigma)+e^{-c/\eps}=KG_*(\tau)+e^{-c/\eps}.
\]
With the bound for $D\notin E$ this gives $|z|_1\le K'(G_*+e^{-c/\eps})$. At
$\tau_\eps$, where $|z|_1=\rho$, it gives $G_*(\tau_\eps)\ge\rho/(2K')$. Afterwards,
as long as $G_*\le\rho$, it gives $|z|_1\le K'(\rho+e^{-c/\eps})<2K'\rho$, which
closes the bootstrap, and $G_*$ reaches $\rho$ within fast time
$2\log(2K')/s_0$. The windows above then extend to $\tau^*_\eps$, which
gives~(c).
\end{proof}

\subsection{Numerical check}\label{sec:eenum}

In seven networks with a two-species siphon, the simulated exit points approach
the predicted tropical values as $\eps$ decreases, also where the Perron integral or the
diagonal block is wrong (Table~\ref{tab:num}). The script \nolinkurl{julia/models/p3_entry_exit.jl} (output
\nolinkurl{julia/results/p3_entry_exit.txt}, about 20 seconds) uses networks with
species $\ce X_1,\ce X_2$ (the siphon) and $\ce U$ (slow, inflow $0\to\ce U$ of
rate $1$, so $u_0(\tau)=\tau$). It integrates the full nonlinear system in slow
time, in the variables $\eps\log x_i$ for the siphon species, whose concentrations
reach $e^{-c/\eps}$, below floating-point range. All simulations in this paper use the
Rosenbrock method Rodas5P of the Julia package OrdinaryDiffEq.jl, with relative and
absolute tolerances $10^{-9}$ and $10^{-11}$ unless stated otherwise. The script
starts at $u=0$ with $x_1=x_2=10^{-2}$, and records $u$, by event detection, when
$\max(x_1,x_2)$ returns to $10^{-2}$; this
exit section differs from $\{|z|_1=\rho\}$ only by $O(\eps)$ in the exit point.
The predictions are computed independently from the network data (block Perron
roots, siphon graph, feeds) with step $10^{-5}$. The script computes what we expect to be the least
solution of \eqref{eq:trop}, as discussed after Conjecture~\ref{conj:cyclic}: it merges the components on a cycle of slow edges,
which then share one $\psi$, and solves the feeds within each step from below, which
removes the spurious solutions through hyperedges. All networks contain
$\ce U+\ce X_1\to\ce U+2\ce X_1$ (rate $1$) and $\ce X_1\to0$ (rate $2$), so
$s_{\{\ce X_1\}}=u-2$, together with $2\ce X_1\to\ce X_1$ and $2\ce X_2\to\ce X_2$ for
saturation.

\begin{table}[ht]
\centering
\small
\begin{tabular}{@{}l>{\raggedright\arraybackslash}p{5.2cm}lllll@{}}
\toprule
& further reactions & \multicolumn{3}{l}{predicted} & \multicolumn{2}{l}{simulated}\\
\cmidrule(lr){3-5}\cmidrule(l){6-7}
&& Perron & blocks & tropical & $\eps=10^{-3}$ & $\eps=10^{-4}$\\
\midrule
I  & $\ce X_1\to\ce X_1+\ce X_2$ $(\tfrac12)$, $\ce X_2\to\ce X_1+\ce X_2$, $\ce X_2\to0$ & 2.71742 & -- & 2.71742 & 2.71695 & 2.71737\\
R1 & $\ce X_2\to\ce X_1+\ce X_2$, $\ce X_2\to0$ $(\tfrac12)$ & 3.32288 & 4 & 3.32288 & 3.31959 & 3.32246\\
R2 & $\ce X_1\to\ce X_1+\ce X_2$, $\ce X_2\to0$ $(\tfrac12)$ & 3.32288 & 4 & 4 & 4.00001 & 4.00000\\
R3 & $\ce X_1\to\ce X_1+\ce X_2$, $\ce X_2\to0$ $(\tfrac1{10})$, $2\ce X_2\to\ce X_1+2\ce X_2$ & 2.62450 & 4 & 2.87178 & 2.87120 & 2.87159\\
R4 & R2 and the slow reaction $\ce X_2\to\ce X_1$ & 3.32288 & 4 & 3.32288 & 3.31958 & 3.32246\\
R5 & R1 with $\ce X_2\to\ce X_1+\ce X_2$ slow & 3.32288 & 4 & 3.32288 & 3.32481 & 3.32316\\
R6 & R3 without $\ce X_1\to\ce X_1+\ce X_2$ & 2.62450 & 4 & 2.87178 & 2.87228 & 2.87170\\
\bottomrule
\end{tabular}
\caption{Exit point $u$ of trajectories entering at $u=0$; rate constants $1$ unless
given in parentheses. ``Blocks'': first zero of $\int s_{\{\ce X_1\}}$, the latest
possible exit (Corollary~\ref{cor:bounds}, proved for acyclic feeding; R3 and R4 are
cyclic); a dash for I, whose only block is the
whole siphon. The predicted values are exact: R1,
R2, R4, R5 Perron $2+\sqrt{7}/2$; R3, R6 tropical $2+\sqrt{0.76}$; R3, R6 Perron
$2+\sqrt{0.39}$.}
\label{tab:num}
\end{table}

\begin{figure}[t]
\centering
\includegraphics[width=12.1cm]{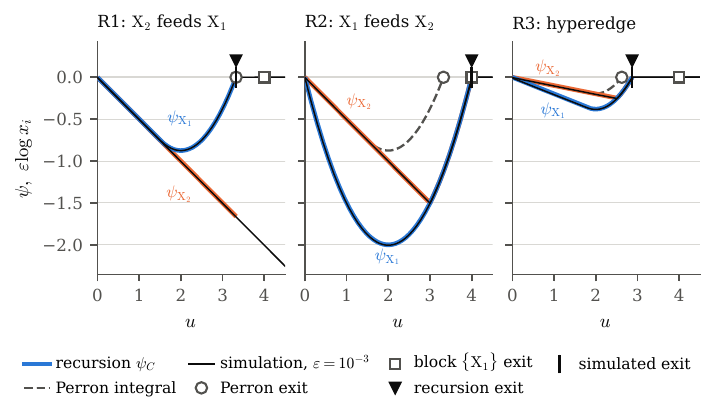}
\caption{Entry--exit for the networks R1, R2 and R3 of Table~\ref{tab:num}, with entry at
$u=0$, on the logarithmic scale $\psi=\eps\log x$. Thick lines: the levels
$\psi_{\ce X_1}$ (blue) and $\psi_{\ce X_2}$ (orange) of the tropical
recursion~\eqref{eq:trop}, up to its exit. Dashed: the Perron integral $\int_0^u s$, up to
its zero; where it coincides with a level it is hidden. Thin black: $\eps\log x_1$ and
$\eps\log x_2$ of the full system for $\eps=10^{-3}$, which lie within $0.01$ of the levels
up to the exit. Markers on $\psi=0$: the exits given by the Perron integral (circle), by
the block $\{\ce X_1\}$ alone (square) and by the recursion (triangle), and the simulated
exit (bar). In R2 and R3 the Perron integral reaches $0$ before the exit.}
\label{fig:tropical}
\end{figure}

\begin{itemize}
  \item \textbf{I} (irreducible). Theorem~\ref{thm:perron}: the exit follows the
        Perron root, error $4.7\cdot10^{-4}$ and $5\cdot10^{-5}$ at $\eps=10^{-3}$
        and $10^{-4}$, i.e.\ $O(\eps)$. The diagonal entry alone ($u-2$) would
        give $4$.
  \item \textbf{R1} ($\ce X_2$ feeds $\ce X_1$). The fed block $\{\ce X_1\}$ changes
        sign. Before $u=\tfrac32$ it decays faster than its feeder, which holds it
        at the level $\psi_{\ce X_2}=-u/2$. After that it grows on its own. The
        tropical recursion and the Perron integral agree, and the diagonal block
        is wrong.
  \item \textbf{R2} ($\ce X_1$ feeds $\ce X_2$). The feeding block changes sign. The
        slowly decaying $\ce X_2$ cannot pass its higher level back to $\ce X_1$, so
        the exit is that of $\ce X_1$ alone. Here the Perron integral is wrong: it
        predicts an exit at $3.32$, where the concentrations are still of order
        $e^{-1.125/\eps}$.
  \item \textbf{R3} (hyperedge). $2\ce X_2\to\ce X_1+2\ce X_2$ feeds $\ce X_1$ with
        weight $2\psi_{\ce X_2}=-0.2u$. This holds $\ce X_1$ up until its own rate
        exceeds $-0.2$, at $u=1.8$, and gives the exit $2+\sqrt{0.76}$. Neither the
        Perron integral nor the linear part of the siphon graph sees this. The
        feeding is cyclic ($\ce X_1\Rightarrow\ce X_2\Rightarrow\ce X_1$), so this
        case tests Conjecture~\ref{conj:cyclic}.
  \item \textbf{R4} (slow back edge). The slow reaction $\ce X_2\to\ce X_1$, of rate
        $\eps$ in fast time, changes the exit of R2 from $4$ to $3.32$. In the limit
        the two components share one $\psi$. As in R3, this is a cyclic case of the
        conjecture.
  \item \textbf{R5, R6} (acyclic tests of Theorem~\ref{thm:tropical}). In R5 the
        feed of R1 is slow, of rate $\eps$; it still changes the exit from the
        block value $4$ to $3.32$. In R6, $\ce X_2$ feeds $\ce X_1$ only through the
        hyperedge, which holds $\ce X_1$ at $2\psi_{\ce X_2}$ as in R3; the feeding
        is acyclic, and the exit is $2+\sqrt{0.76}$, not the Perron value.
\end{itemize}
Figure~\ref{fig:tropical} shows R1--R3 on the logarithmic scale: the simulated
$\eps\log x_i$ follow the levels $\psi_C$, and in R2 and R3 the Perron integral
reaches $0$ before the exit.

\paragraph{Olsen's model.} The recursion can also be read off in Olsen's model
\citep{Olsen1983}, although it lies outside our setting in two ways. It has an inflow $0\to\ce X$ of rate $k_6>0$, so $\{\ce X,\ce Y\}$ is not a
siphon of the whole network and (S1) fails. In the scaling of
\citet{KuehnSzmolyan2015} this inflow is a breaking reaction in the sense of
Sections~\ref{sec:break} and~\ref{sec:floor}, and their Cor.~4.5 concerns the
regime in which it is exponentially small; there the floor of
Section~\ref{sec:floor} lies below the passage. And the fast reaction
$\ce A+\ce B+\ce Y\to3\ce X$ changes $\ce A$ and $\ce B$, so (S2) fails
(Remark~\ref{rem:landing}). The siphon
$\{\ce X,\ce Y\}$ has the fast edge $\ce Y\to\ce X$ (through
$\ce A+\ce B+\ce Y\to3\ce X$) and the hyperedge $2\ce X\to2\ce Y$, so the feeding is
cyclic. Along the passages computed by \citet{KuehnSzmolyan2015}, the rate
$s_{\{\ce Y\}}=-(k_3ab+k_5)$ of the $\{\ce Y\}$ block stays below the rate
$s_{\{\ce X\}}=k_1b-k_4$ of the $\{\ce X\}$ block. So $\psi_{\ce X}$ follows its own
block, $\psi_{\ce Y}$ soon equals $2\psi_{\ce X}<\psi_{\ce X}$, and the exit is the
one used in \citet[Cor.~4.5]{KuehnSzmolyan2015}. The Perron integral, the tropical
recursion and the $\{\ce X\}$ block all agree there.

%======================================================================
\section{Relaxation oscillations through the face}\label{sec:ro}
%======================================================================

We keep (S1) and (S2). In addition, the fast species are exactly the species of
$\Sigma$, and there is one slow species $\ce U$. So
\[
  \eps z'=M(u,z,\eps)z,\qquad u'=f(u,z,\eps),\qquad z\in\R^k_{\ge0},\ u>0 .
\]
The critical manifold contains the face $\{z=0\}$ and the interior set
$\{M(u,z,0)z=0,\ z>0\}$. It may contain further pieces on the faces of smaller
siphons; the cycles below do not use them. On the face the transverse block is $\Lambda(u)$, and the
slow speed is $f(u,0,0)$, given by the slow reactions with no siphon species among
their reactants.

\begin{definition}[Fold--face cycle]\label{def:ff}
A \emph{fold--face cycle} is a closed singular orbit $\Gamma$ made of four pieces,
with the following properties. We take the face flow increasing $u$; the other
orientation is symmetric.
\begin{enumerate}[label=(C\arabic*),nosep]
  \item \emph{Branch.} The interior critical set contains a curve $z=h(u)>0$ for
        $u\in[u_F,u_1]$. It is normally attracting for $u>u_F$: the fast Jacobian
        $J=D_z(M z)$ is Hurwitz. It ends at a regular fold at $u_F$: there $J$ has
        a simple zero eigenvalue, with right and left null vectors $r_F$ and
        $\ell_F$, the other eigenvalues have negative real part,
        $\ell_F\cdot D_z^2(Mz)(r_F,r_F)\ne0$ and $\ell_F\cdot\partial_u(Mz)\ne0$,
        and the slow flow on the branch points towards the fold
        ($f(u,h(u),0)<0$, including at $u_F$).
  \item \emph{Drop.} In the layer problem at $u=u_F$, the fast fibre that leaves
        the fold point along its centre direction (the \emph{drop fibre})
        converges to the face point $(0,u_F)$, inside the open orthant. The face
        is attracting there: $s(\Lambda(u_F))<0$.
  \item \emph{Face.} The feeding relation is acyclic. The face speed is positive on
        $[u_F,u_{\rm out}]$, where $u_{\rm out}=u_0(\tau_{\rm ex})$ is the exit
        point of Theorem~\ref{thm:tropical} for entry at $u_F$ with all entry
        values $a_C=0$ (the Perron exit of Theorem~\ref{thm:perron} if $\Lambda$ is
        irreducible). The exit is transversal (Definition~\ref{def:transversal}),
        with exiting components $E$ and self-exiting components $E^*$.
  \item \emph{Rise.} There is $\rho_1>0$ such that every layer trajectory at
        $u=u_{\rm out}$ from a point $z\ge0$ with $|z|_1\le\rho_1$, $z_D=0$ for
        $D\notin E$, and $z_C\ne0$ for some $C\in E^*$, converges to
        $h(u_{\rm out})$. Moreover $u_F<u_{\rm out}<u_1$.
\end{enumerate}
Its pieces are the branch from $u_{\rm out}$ down to $u_F$, the drop, the face
segment $[u_F,u_{\rm out}]$, and the rise.
\end{definition}

\begin{theorem}[Relaxation oscillation through a siphon face]\label{thm:ro}
Let $\Gamma$ be a fold--face cycle. For all small $\eps>0$ there is a unique periodic
orbit $\Gamma_\eps$ in a neighbourhood of $\Gamma$. It is orbitally asymptotically
stable, with nontrivial Floquet multipliers $O(e^{-c/\eps})$, and $\Gamma_\eps\to\Gamma$
in the Hausdorff distance. Its minimum and maximum of $u$ tend to $u_F$ and
$u_{\rm out}$, and its period tends to
\[
  T_0=\int_{u_F}^{u_{\rm out}}\frac{du}{f(u,0,0)}+\int_{u_F}^{u_{\rm out}}\frac{du}{|f(u,h(u),0)|}.
\]
\end{theorem}

The proof uses one estimate beyond Section~\ref{sec:ee}: the passage along the face
is $C^1$, with derivative $O(1/\eps)$. Let $\varpi=\one^{\top}(-\Lambda(u_F))^{-1}>0$,
and
\[
  G_*(z,u)=\sum_{C\in E^*}\ell_C(u)\cdot z_C ,
\]
where $\ell_C(u)$ is the left Perron vector of $\Lambda_{CC}(u)$ with
$\ell_C\cdot\one=1$.
\begin{itemize}[nosep]
  \item Let $\rho_e=\rho\min_i\varpi_i$, so that $\varpi\cdot z=\rho_e$ implies
        $|z|_1\le\rho$. Let $p_\rho=(z^\rho,u_F)$ be the first point of the drop
        fibre with $\varpi\cdot z=\rho_e$. By (C2), $z^\rho>0$. The fibre crosses
        $\{\varpi\cdot z=\rho_e\}$ transversally: along it, in fast time,
        $(\varpi\cdot z)^{\boldsymbol\cdot}=\varpi\Lambda(u_F)z+O(|z|_1^2)=-|z|_1+O(|z|_1^2)<0$.
  \item For a small $\gamma_0>0$, let $\mathcal K_E$ be the set of points of
        $\{\varpi\cdot z=\rho_e\}$ with $|z-z^\rho|\le\gamma_0$ and
        $|u-u_F|\le\gamma_0$. It is a compact subset of the open orthant.
  \item Let $L$ bound $|du_{\rm out}/du_{\rm in}|$ near $u_F$, and
        $\mathcal X=\{G_*=\rho,\ |u-u_{\rm out}|\le2L\gamma_0\}$.
\end{itemize}

\begin{lemma}[Face passage]\label{lem:face}
Assume (C2) and (C3). For small $\rho$, $\gamma_0$ and $\eps$, define the passage
map $P_\eps:\mathcal K_E\to\mathcal X$ by following the trajectory to the first time
$\tau^*_\eps\ge\tau_0$ at which $G_*=\rho$, for a fixed small $\tau_0$. Then $P_\eps$
is well defined and $C^1$, and:
\begin{enumerate}[label=(\alph*),nosep]
  \item at $P_\eps(p)$, $u=u_{\rm out}(u(p))+o(1)$ uniformly on $\mathcal K_E$,
        $|z|_1\le K\rho$, and $|z_D|\le e^{-c/\eps}$ for every component
        $D\notin E$;
  \item $\|DP_\eps\|\le K/\eps$ on $\mathcal K_E$.
\end{enumerate}
\end{lemma}

\begin{proof}
Lemmas~\ref{lem:apriori} and~\ref{lem:integral} apply, with $K_0$ the
$z$-projection of $\mathcal K_E$ and $U_0=[u_F-\gamma_0,u_F+\gamma_0]$, with $E$
and $E^*$ those of the exit from $u_F$. Part~(a) follows
from Theorem~\ref{thm:tropical}, from Lemma~\ref{lem:integral}(a) and~(c), and from
the bound on the components $D\notin E$ in the proof of Lemma~\ref{lem:integral}.

For (b), work in the coordinates $(\xi,u)$, $\xi_i=\log z_i$. In them
$\eps\xi_i'=\Xi_i(\xi,u,\eps)$, a sum of terms $\Theta=n\kappa(u,\eps)e^{y\cdot\xi-\xi_i}$,
and $u'=f(u,e^\xi,\eps)$.

\emph{Weights.} Define $\nu_C\ge1$ along the acyclic order by
\[
  \nu_C=\max\bigl(1,\ \max\{y_r\cdot\nu:\ r\text{ a feed of }C\}\bigr),
\]
where $y\cdot\nu=\sum_jy_j\nu_{c(j)}$. Put $\nu_i=\nu_{c(i)}$ and
$\|x\|_\nu=\max_i|x_i|/\nu_i$. The associated logarithmic norm is
\[
  \mu_\nu(J)=\max_i\Bigl(J_{ii}+\sum_{j\ne i}|J_{ij}|\,\nu_j/\nu_i\Bigr),
\]
and solutions of $x'=J(\tau)x$ satisfy $D^+\|x\|_\nu\le\mu_\nu(J)\|x\|_\nu$.

\emph{The $\xi$-Jacobian.} A term $\Theta$ of $\Xi_i$ has
$\partial_{\xi}\Theta=\Theta\,(y-e_i)$.
\begin{itemize}[nosep]
  \item A production term ($n>0$) adds $\Theta\,(y\cdot\nu/\nu_i-1)$ to row $i$ of
        $\mu_\nu$. This is $\le0$ when $y\cdot\nu\le\nu_i$. That covers $y=e_i$, the
        linear edges inside a component, and, by the choice of $\nu$, all feeds
        from other components.
  \item A consumption term ($n<0$, $y_i\ge1$) adds at most
        $|\Theta|\sum_{j\ne i}y_j\nu_j/\nu_i$, which is $0$ when $y=e_i$.
\end{itemize}
The remaining terms are of two kinds:
\begin{itemize}[nosep]
  \item production terms with $|y|\ge2$ and a reactant $X_j$ in the component of
        $X_i$, equal to $n\kappa z^{y-e_j}(z_j/z_i)$;
  \item consumption terms with $|y|\ge2$, equal to $-|n|\kappa z^{y-e_i}$.
\end{itemize}
By Lemma~\ref{lem:apriori}(a), both are $O(|z|_1)$. Hence
$\mu_\nu(D_\xi\Xi)\le K|z|_1$.

\emph{Other derivatives.}
\begin{itemize}[nosep]
  \item $\partial_u\Xi_i=\sum(\partial_u\kappa/\kappa)\,\Theta$ is bounded, by
        Lemma~\ref{lem:apriori}(b).
  \item $\partial_\xi f=(\partial_zf)\operatorname{diag}(z)=O(|z|_1)$, and
        $\partial_uf$ is bounded.
\end{itemize}

\emph{Gronwall.} Let $(\delta\xi,\delta u)$ solve the variational equation along
a trajectory from $\mathcal K_E$, and set $A=\eps\|\delta\xi\|_\nu$, $b=|\delta u|$. Then
\[
  D^+A\le\frac{K|z|_1}{\eps}A+Kb,\qquad D^+b\le Kb+\frac{K|z|_1}{\eps}A .
\]
By Lemma~\ref{lem:integral}(c), $\int_0^{\tau^*_\eps}|z|_1/\eps\le K$. So
$A+b\le K\,(A(0)+b(0))$ on $[0,\tau^*_\eps]$. That is,
\[
  \|\delta\xi\|\le K\bigl(\|\delta\xi(0)\|+|\delta u(0)|/\eps\bigr),\qquad
  |\delta u|\le K\bigl(\eps\|\delta\xi(0)\|+|\delta u(0)|\bigr).
\]

\emph{Back to $z$, and the sections.}
\begin{itemize}[nosep]
  \item On $\mathcal K_E$, $z_i\ge c>0$, so $|\delta\xi(0)|\le K|\delta z(0)|$. At
        $\tau^*_\eps$, $|\delta z|\le|z|\,|\delta\xi|\le K\rho|\delta\xi|$. So the
        derivative of the flow up to time $\tau^*_\eps$ is $O(1/\eps)$.
  \item The exit time $\tau^*_\eps$ is a $C^1$ function of $p\in\mathcal K_E$, by
        the implicit function theorem, provided the crossing of $\{G_*=\rho\}$ is
        transversal. Since feeds are nonnegative,
        $\ell_C\cdot(\eps z_C')\ge(s_C-K|z|_1-K\eps)Z_C$, and $s_C\ge s_0>0$ for
        $C\in E^*$ near the exit. Hence $\frac{d}{d\tau}G_*\ge s_0\rho/(4\eps)$ at
        $\tau^*_\eps$.
  \item The vector field is $O(1/\eps)$. So the exit time has bounded
        derivative, and $\|DP_\eps\|\le K/\eps$. \qedhere
\end{itemize}
\end{proof}

\begin{proof}[Proof of Theorem~\ref{thm:ro}]
Fix $u_2<u_3$ in $(u_F,u_{\rm out})$, and let $\mathcal V$ be a neighbourhood of
the branch over $[u_F+\gamma_1,u_{\rm out}+\gamma_1]$, with $u_{\rm out}+\gamma_1<u_1$,
on which Fenichel's theorem applies \citep{Fenichel1979}. In Fenichel coordinates
$(\chi,\bar u)\in\R^k\times\R$ on $\mathcal V$:
\begin{itemize}[nosep]
  \item the slow manifold is $S_\eps=\{\chi=0\}$;
  \item the base point of the stable fibre moves by the reduced flow,
        $\bar u'=\bar g(\bar u,\eps)$;
  \item $\eps\chi'=\tilde A(\chi,\bar u,\eps)\chi$, where $\tilde A$ has spectrum in
        $\{\operatorname{Re}\le-2c\}$ \citep[Ch.~3]{Kuehn2015}.
\end{itemize}
Let $\Delta$ and $\Delta_3$ be the sections $\{\bar u=u_2\}$ and $\{\bar u=u_3\}$ of
$\mathcal V$, and
let $B\subset\Delta$ be a closed ball of fixed radius around
$p_0=\Gamma\cap\Delta$. Write the return map as
$\Pi=\Pi_4\circ\Pi_3\circ\Pi_2\circ\Pi_1$. The four maps are as follows.

\emph{$\Pi_1$: from $B$, past the fold, to a section $\Delta_f$ transverse to the
drop fibre.}
\begin{itemize}[nosep]
  \item Along $\mathcal V$ the fibres contract, and the base point reaches the fold
        region in slow time $O(1)$.
  \item At the fold, the fast eigenvalues other than the zero one have negative
        real part. So a centre manifold of dimension two in $(z,u)$, plus $\eps$,
        attracts the flow at a fast exponential rate. On it the flow is a planar
        slow--fast system with a fold, nondegenerate by the conditions in (C1).
        The slow flow points into the fold, so the fold theorem of
        \citet{KrupaSzmolyan2001fold} applies.
  \item The transition map is defined on $B$, has bounded derivative (it is in
        fact a contraction), and maps $B$ into an $O(\eps^{2/3})$-neighbourhood of
        the drop fibre.
\end{itemize}

\emph{$\Pi_2$: from $\Delta_f$ to $\mathcal K_E$.} The drop fibre reaches
$\{\varpi\cdot z=\rho_e\}$ in finite fast time, transversally. So $\Pi_2$ is a regular
perturbation with bounded derivative, and its image lies within $O(\eps^{2/3})$ of
$p_\rho$, hence in $\mathcal K_E$.

\emph{$\Pi_3=P_\eps$}, by Lemma~\ref{lem:face}.

\emph{$\Pi_4$: from a neighbourhood of $\mathcal X$ to $\Delta$.}
\begin{itemize}[nosep]
  \item By Lemma~\ref{lem:face}(a), the image $P_\eps(\Pi_2(\Pi_1(B)))$, which lies
        within $O(\eps^{2/3})$ of the $u$-level $u_{\rm out}(u_F)$ at entry, lies
        within $o(1)$ of the compact set
        $\mathcal X_0=\{G_*=\rho,\ u=u_{\rm out},\ z_D=0\text{ for }D\notin E,\ |z|_1\le K\rho,\ z\ge0\}$.
  \item On $\mathcal X_0$ some $z_C$, $C\in E^*$, is nonzero, since $G_*=\rho$. So
        for $K\rho\le\rho_1$, by (C4), every layer trajectory from $\mathcal X_0$
        converges to $h(u_{\rm out})$. Since $h(u_{\rm out})$ is hyperbolic and
        $\mathcal X_0$ is compact, all such trajectories enter a compact subset of
        $\mathcal V$ within a fixed fast time $t_1$, and so do those from a
        neighbourhood of $\mathcal X_0$. Over that time the flow has bounded
        derivative.
  \item Inside $\mathcal V$, $\bar u$ moves from $u_{\rm out}+o(1)$ to $u_3$ and then to $u_2$
        in slow time at least $t_2>0$. The time to reach $\Delta_3$ depends only on
        the base point, through the reduced flow.
  \item So the $\chi$-coordinate on arrival in $\Delta_3$ has derivative
        $O(e^{-2ct_2/\eps}/\eps)$ with respect to the initial point. The map
        $\Delta_3\to\Delta$ along $\mathcal V$ has bounded derivative.
\end{itemize}
Hence $\|D\Pi_4\|\le Ce^{-ct_2/\eps}$. Moreover $\Pi_4$ maps into an
$O(e^{-ct_2/\eps})$-neighbourhood of $p_\eps=S_\eps\cap\Delta$, which lies in $B$.

\emph{Fixed point.} $\Pi$ maps $B$ into itself, and on $B$
\[
  \|D\Pi\|\le C\cdot C\cdot\frac C\eps\cdot Ce^{-ct_2/\eps}\le e^{-c'/\eps}<1 .
\]
$B$ is convex, so $\Pi$ is a contraction on $B$, and it has a unique fixed point
$p_\eps^*$. Let $\Gamma_\eps$ be its orbit.
\begin{itemize}[nosep]
  \item \emph{Uniqueness.} A periodic orbit in a thin enough tubular neighbourhood
        of $\Gamma$ crosses $\Delta$ inside $B$. If it winds $n$ times around the
        neighbourhood, its crossing is a fixed point of $\Pi^n$, which is also a
        contraction. So it is $p^*_\eps$.
  \item \emph{Stability.} The nontrivial Floquet multipliers are the eigenvalues
        of $D\Pi(p^*_\eps)$, of size at most $e^{-c'/\eps}$. So $\Gamma_\eps$ is
        orbitally asymptotically stable.
\end{itemize}

\emph{Limits.} Piece by piece:
\begin{itemize}[nosep]
  \item $\Gamma_\eps$ is $O(\eps)$-close to the branch inside $\mathcal V$;
  \item it is $O(\eps^{1/3})$-close to the fold point and the drop fibre between
        the branch and $\Delta_f$ \citep{KrupaSzmolyan2001fold};
  \item between $\Delta_f$ and $\mathcal K_E$ it is close to the drop fibre;
  \item on the face, $|z|\le K\rho$ and $|u(\tau)-u_0(\tau)|=O(\eps)$, by Gronwall
        and Lemma~\ref{lem:integral}(b);
  \item between $\mathcal X$ and $\mathcal V$ it is close to a layer trajectory from
        $\mathcal X_0$.
\end{itemize}
So $\limsup_{\eps\to0}d_H(\Gamma_\eps,\Gamma)\le\omega(\rho)$, with
$\omega(\rho)\to0$ as $\rho\to0$. The orbit $\Gamma_\eps$ does not depend on
$\rho$, by uniqueness, so $d_H(\Gamma_\eps,\Gamma)\to0$. The extreme values of
$u$ follow, the minimum $u_F+O(\eps^{2/3})$ coming from the fold passage.

\emph{Period.} The time spent in $\mathcal V$ tends to the branch integral, since
$\bar u'=\bar g(\bar u,\eps)$ is the reduced flow up to $O(\eps)$. The time between
$\Delta_f$ and $\mathcal K_E$, and between $\mathcal X$ and $\mathcal V$, is
$O(\eps)$, and the fold passage takes slow time
$o(1)$. The face time is $\tau^*_\eps\to\tau_{\rm ex}=\int_{u_F}^{u_{\rm out}}du/f(u,0,0)$.
Letting $\rho,\gamma_0,\gamma_1\to0$ gives $T_\eps\to T_0$.
\end{proof}

\begin{remark}\label{rem:ro}
\begin{itemize}[nosep]
  \item For $k=1$ in the plane this is the type of theorem proved by
        \citet{AiSadhu2020}, and for specific predator--prey models by
        \citet{LiLiWang2023}. The new parts are coupled siphon species ($k\ge2$) and the reading of the
        exit point from the network.
  \item In the settings of \citet{Hsu2019} and \citet{HsuRuan2021} the entry point
        moves along a family of landing points. A cycle is then selected by an
        integral condition, and its stability needs a Floquet-type condition.
        Here the fold fixes the entry, so both come for free.
  \item Exiting components that are not self-exiting, such as $\ce X_2$ in
        $\Net3$ below, are carried to the exit by their feeders and decay on their
        own there. The exit section $\{G_*=\rho\}$ counts only the self-exiting
        components, which grow; (C4) asks the layer flow to rise from points where
        one of them is present.
  \item With cyclic feeding, (C3) is not available: the exit point rests on
        Conjecture~\ref{conj:cyclic}, and the weights in the proof of
        Lemma~\ref{lem:face} do not exist, since a cycle through a hyperedge needs
        $\nu_C\ge2\nu_C$. Network~$\Net2$ below is such a case. The numerics
        support the conclusion there too.
\end{itemize}
\end{remark}

\subsection{Numerical check}\label{sec:ronum}

The script \nolinkurl{julia/models/p3_relaxation.jl} (output
\nolinkurl{julia/results/p3_relaxation.txt}, about 30 seconds) integrates the full
networks in slow time. It uses the variables $\eps\log x$ for the siphon species, which
fall to $e^{-c/\eps}$ on the face, and reports the range of the slow species and the
period over the last three cycles. The exit predictions come from the same code as
in Section~\ref{sec:eenum}. Table~\ref{tab:ro} lists the results.

Network $\Net1$, with fast $\ce X$ and slow $\ce Y$, is the planar example of
Section~\ref{sec:planar}.
Networks $\Net2$ and $\Net3$ have the siphon $\{\ce X_1,\ce X_2\}$ and the slow
species $\ce U$. Their fast reactions are $\ce U+\ce X_1\to\ce U+2\ce X_1$ $(1)$,
$\ce X_1\to0$ $(4)$, $2\ce X_1\to3\ce X_1$ $(3)$, $3\ce X_1\to2\ce X_1$ $(1)$,
$\ce X_1\to\ce X_1+\ce X_2$ $(1)$, $\ce X_2\to0$ $(\tfrac1{10})$ and
$2\ce X_2\to\ce X_2$ $(1)$. Network $\Net2$ adds the hyperedge
$2\ce X_2\to\ce X_1+2\ce X_2$ $(1)$. The slow reactions are $0\to\ce U$ $(1)$ and
$\ce U+\ce X_1\to\ce X_1$ $(2)$.
\begin{itemize}[nosep]
  \item On the face, $s_{\{\ce X_1\}}=u-4$ changes sign, while $\ce X_2$ decays
        slowly at rate $\tfrac1{10}$.
  \item Network $\Net3$ has acyclic feeding ($\ce X_1$ feeds $\ce X_2$), and its
        exit is transversal. Both components exit, but only $\{\ce X_1\}$ is
        self-exiting, with $s_{\{\ce X_1\}}=u_{\rm out}-4>0$, while
        $s_{\{\ce X_2\}}=-\tfrac1{10}$. It satisfies the hypotheses of
        Theorem~\ref{thm:ro}, checked below, so the theorem applies. The exit is that of the block
        $\{\ce X_1\}$, and the Perron integral is wrong (as in case R2 of
        Section~\ref{sec:eenum}). Since $\ce X_2$ acts on neither $\ce X_1$ nor
        $\ce U$, the cycle is that of the planar system in $(x_1,u)$; the
        multi-species Perron integral misses it.
  \item In $\Net2$ the hyperedge holds $\ce X_1$ at $2\psi_{\ce X_2}$, and the
        feeding is cyclic. The exit differs from both the Perron integral and the
        block.
\end{itemize}

\paragraph{Hypotheses for $\Net3$.} Here $\eps x_1'=x_1\varphi$ with
$\varphi=u-4+3x_1-x_1^2$, $\eps x_2'=x_1-\tfrac1{10}x_2-x_2^2$ and $u'=1-2ux_1$.
(S1) and (S2) hold: every reaction that produces $\ce X_1$ or $\ce X_2$ has one of
them among its reactants, and no fast reaction changes $\ce U$. The interior
critical set is $u=H(x_1)=4-3x_1+x_1^2$, $x_2=h_2(x_1)=(\sqrt{1+400x_1}-1)/20$.
\begin{itemize}[nosep]
  \item (C1) The fast Jacobian on it is lower triangular, with diagonal
        $x_1(3-2x_1)$ and $-\tfrac1{10}-2x_2<0$. So the branch $x_1>\tfrac32$, on
        which $H$ increases, is attracting. It ends at the fold $x_1=\tfrac32$,
        $u_F=\tfrac74$, where the zero eigenvalue is simple, with
        $\ell_F=(1,0)$ and $r_F=(a,1)$, $a=\tfrac1{10}+2x_2$. There
        $\ell_F\cdot D_z^2(Mz)(r_F,r_F)=-3a^2\ne0$ and
        $\ell_F\cdot\partial_u(Mz)=x_1=\tfrac32\ne0$. On the branch $x_1\ge\tfrac32$
        and $u\ge\tfrac74$, so $u'\le1-2\cdot\tfrac74\cdot\tfrac32<0$.
  \item (C2) At $u=u_F$, $\varphi=-(x_1-\tfrac32)^2<0$ for $x_1\ne\tfrac32$. So the
        drop fibre takes $x_1$ to $0$, and then $x_2\to0$. The face is attracting
        there: $\Lambda(u_F)$ is lower triangular with diagonal $-\tfrac94$ and
        $-\tfrac1{10}$.
  \item (C3) The feeding is acyclic and the face speed is $u'=1$. The exit at
        $u_{\rm out}=\tfrac{25}4$ is transversal, as stated above.
  \item (C4) At $u=u_{\rm out}$, $\varphi=-(x_1^2-3x_1-\tfrac94)$ is positive on
        $(0,x_1^+)$ and negative beyond, with $x_1^+=\tfrac32(1+\sqrt2)$. So every
        layer trajectory with $x_1>0$ has $x_1\to x_1^+$, and then
        $x_2\to h_2(x_1^+)$, the unique positive equilibrium of the
        $x_2$-equation, which attracts. This is the branch point over $u_{\rm out}$,
        and $u_F<u_{\rm out}<u_1$ for any $u_1>\tfrac{25}4$.
\end{itemize}

\begin{table}[ht]
\centering
\small
\begin{tabular}{@{}lllll@{}}
\toprule
network & & lower turning point & upper turning point & period\\
\midrule
$\Net1$ (planar, $y$) & singular cycle & $y_{\rm out}=0.40638$ & $y_F=2$ & $4.3310$\\
 & $\eps=10^{-2}$ & $0.38243$ & $2.08339$ & $4.6887$\\
 & $\eps=10^{-3}$ & $0.40004$ & $2.01987$ & $4.4160$\\
 & $\eps=10^{-4}$ & $0.40485$ & $2.00454$ & $4.3505$\\
\midrule
$\Net2$ (hyperedge, $u$) & singular cycle & $u_F=0.82823$ & $u_{\rm out}=5.10847$ & $4.7067$\\
 & \quad Perron / block & & $4.79016$ / $7.17177$ & \\
 & $\eps=10^{-2}$ & $0.72925$ & $5.10299$ & $4.9145$\\
 & $\eps=10^{-3}$ & $0.80308$ & $5.10966$ & $4.7549$\\
 & $\eps=10^{-4}$ & $0.82250$ & $5.10906$ & $4.7175$\\
\midrule
$\Net3$ (acyclic, $u$) & singular cycle & $u_F=1.75$ & $u_{\rm out}=6.25$ & $4.7634$\\
 & \quad Perron / block & & $4.66332$ / $6.25$ & \\
 & $\eps=10^{-2}$ & $1.57431$ & $6.41881$ & $5.1812$\\
 & $\eps=10^{-3}$ & $1.70743$ & $6.29191$ & $4.8621$\\
 & $\eps=10^{-4}$ & $1.74035$ & $6.25932$ & $4.7855$\\
\bottomrule
\end{tabular}
\caption{Relaxation oscillations through a siphon face. In $\Net1$ the face is
traversed with $y$ decreasing, so its exit point is the lower turning point; in
$\Net2$ and $\Net3$ it is the upper one. The predicted exits are given in closed form: $\Net1$,
the root below $2$ of $\ln(y/2)=y-2$; $\Net2$, $4+\sqrt{0.04+0.4\,(3.8-u_F)}$; $\Net3$,
$4+(4-u_F)$. The fold values $u_F$ are
computed on the interior curve, parametrised by $x_1$. ``Perron / block'': the
exits that the Perron integral and the block $\{\ce X_1\}$ alone would give.
Periods in slow time.}
\label{tab:ro}
\end{table}

\begin{itemize}[nosep]
  \item The turning point at the fold converges at the rate $\eps^{2/3}$. The
        error ratios per decade of $\eps$ are $3.9$--$4.4$, approaching
        $10^{2/3}=4.64$. This error is the delay at the fold.
  \item In $\Net1$ and $\Net3$ the exit point inherits that error, shifted through
        the entry. In $\Net2$ the exit is within $1.2\cdot10^{-3}$ already at
        $\eps=10^{-3}$: there $du_{\rm out}/du_F=-0.18$, against $-1$ in
        $\Net3$, so the delayed entry matters less.
  \item The periods converge to $T_0$ at the same rate.
  \item In $\Net2$ the Perron prediction misses the exit point by $0.32$ and the
        block by $2.06$. In $\Net3$ the Perron prediction misses by $1.59$.
\end{itemize}

%======================================================================
\section{Planar fold--face networks}\label{sec:planar}
%======================================================================

\subsection{Conditions on the reactions}

Let $\Sigma=\{\ce X\}$, with $\ce X$ fast and $\ce Y$ slow, and let the system
have the form $\eps x'=x\,\varphi(x,y)$, $y'=g(x,y)$ in slow time. This is the form
of a network in the setting of Section~\ref{sec:setting}, where (S2) asks that no
fast reaction change $\ce Y$. It is also the form of the species scaling of
Section~\ref{sec:enum}, in which the time scales belong to the species,
and every reaction enters both equations: $\varphi$
collects the reactions that change $\ce X$, and $g$ those that change $\ce Y$. The
proof of Theorem~\ref{thm:ro} uses only this form. Write the reactions that change
$\ce X$ as $a_r\ce X+b_r\ce Y\to\cdots$, with net change $L_r$ of $\ce X$ and rate
constant $k_r$. Those with $L_r>0$ have $a_r\ge1$, because $\{\ce X\}$ is a siphon,
and $L_r<0$ forces $a_r\ge1$ anyway. So
\[
  \varphi(x,y)=\sum_rL_rk_rx^{a_r-1}y^{b_r}.
\]
Assume it is affine in $y$ with negative slope: every such reaction has
$b_r\le1$, and those with $b_r=1$ consume $\ce X$. Then
$\varphi=A(x)+yB(x)$ with $B<0$ on $x\ge0$. Write
\[
  H(x)=-\frac{A(x)}{B(x)},\qquad
  \varphi=B(x)\bigl(y-H(x)\bigr),\qquad
  P=A'B-AB' .
\]
The interior critical set is the graph $y=H(x)$. On it,
$\varphi_x=P/B$ and $H'=-P/B^2$. So the branch is attracting exactly where
$P>0$, which is where $H$ decreases. The face has transverse rate
$\lambda(y)=\varphi(0,y)=A(0)+yB(0)$, and $g(x,y)$ is the slow speed of $\ce Y$.

\begin{proposition}[Planar fold--face cycles]\label{prop:planar}
Assume the following.
\begin{enumerate}[label=(N\arabic*),nosep]
  \item $A(0)>0$ and $B(0)<0$. That is, some reaction $\ce X\to\cdots$ with no
        $\ce Y$ produces $\ce X$, and some reaction $\ce X+\ce Y\to\cdots$ consumes
        $\ce X$.
  \item $P(0)<0$, and $P$ has exactly one positive root $x_F$, which is simple. For
        instance, the coefficients of $P$ have exactly one sign change.
  \item $H(x)<0$ for large $x$. For instance, the net $\ce X$-change of the
        $\ce Y$-free reactions of highest $\ce X$-order is negative.
  \item $g(0,y)<0$ for $y>0$: the reactions without $\ce X$ consume $\ce Y$
        overall.
  \item $g(x,H(x))>0$ for $x\in[x_F,x_+]$, where $x_+>x_F$ solves
        $H(x_+)=y_{\rm out}$ and $y_{\rm out}$ is defined below.
\end{enumerate}
Then $H$ rises from $y_T=H(0)$ to its maximum $y_F=H(x_F)$ and then decreases. The
face is attracting for $y>y_T$ and repelling for $y<y_T$. The equation
\[
  \int_{y_F}^{y_{\rm out}}\frac{\lambda(y)}{g(0,y)}\,dy=0
\]
has a unique solution $y_{\rm out}\in(0,y_T)$, and there is a fold--face cycle with
$u_F=y_F$ and $u_{\rm out}=y_{\rm out}$, run with $y$ decreasing on the face. By
Theorem~\ref{thm:ro}, the network has a unique attracting relaxation oscillation
near it for all small $\eps>0$.
\end{proposition}

\begin{proof}
By (N2), $P<0$ on $(0,x_F)$ and $P>0$ beyond. So $H$ increases up to $x_F$ and
decreases afterwards. By (N1),
$y_T=-A(0)/B(0)>0$, and (N3) makes $H$ negative eventually. The fold at $x_F$ is
regular: $P'(x_F)\ne0$ gives $H''(x_F)\ne0$, and $g(x_F,y_F)>0$ by (N5).

\emph{Drop.} For $0\le x<x_F$, $\varphi(x,y_F)=B(x)(y_F-H(x))<0$. So the fibre
from the fold goes to the face. There
$\lambda(y_F)=B(0)(y_F-y_T)<0$, so the face is attracting.

\emph{Face.} The face is attracting for $y>y_T$ and repelling below. By (N4) the
face flow decreases $y$. In slow time the integral is $\int\lambda\,d\tau$, with
$d\tau=dy/g(0,y)>0$. It is negative while $y>y_T$, where $\lambda<0$, and increases
once $y<y_T$. Consuming $\ce Y$ needs $\ce Y$ among the reactants, so $g(0,y)=O(y)$.
Since $\lambda(0)=A(0)>0$, the integral tends to $+\infty$ as $y\to0$. So it has
exactly one zero $y_{\rm out}\in(0,y_T)$.

\emph{Rise.} For $y_{\rm out}<y_T$, $\varphi(x,y_{\rm out})=B(x)(y_{\rm out}-H(x))>0$
for $0\le x<x_+$, because $H>y_T>y_{\rm out}$ on $(0,x_F]$ and $H$ decreases to
$y_{\rm out}$ at $x_+$. So the fibre rises to $x_+$ on the attracting branch.

By (N5) the slow flow carries the branch from $x_+$ to the fold.
\end{proof}

\begin{corollary}[Predator--prey form]\label{cor:pp}
If, in addition, the slow speed has the form $g=y\,G(x)$ with $G$ increasing (for
example $G=cx-m$, from $\ce X+\ce Y\to\ce X+2\ce Y$ and $\ce Y\to0$), then (N4)
reads $G(0)<0$, and (N5) reads $G(x_F)>0$. That is, the coexistence equilibrium
$x^*=G^{-1}(0)$ lies on the repelling side of the fold, $x^*<x_F$.
\end{corollary}

\begin{remark}\label{rem:explosion}
As $x^*$ increases through $x_F$, the equilibrium crosses the fold and becomes a
canard point in the sense of \citet{BeldimanHalanay2026}. We expect the canard
explosion to connect small Hopf cycles to the fold--face relaxation oscillation,
whose large cycles pass along the face; compare \citet{LiLiWang2023}. The canard
cycles along the repelling middle branch then have to pass the transcritical point,
and we do not prove this here (Section~\ref{sec:disc}).
\end{remark}

\begin{figure}[t]
\centering
\includegraphics[width=12.1cm]{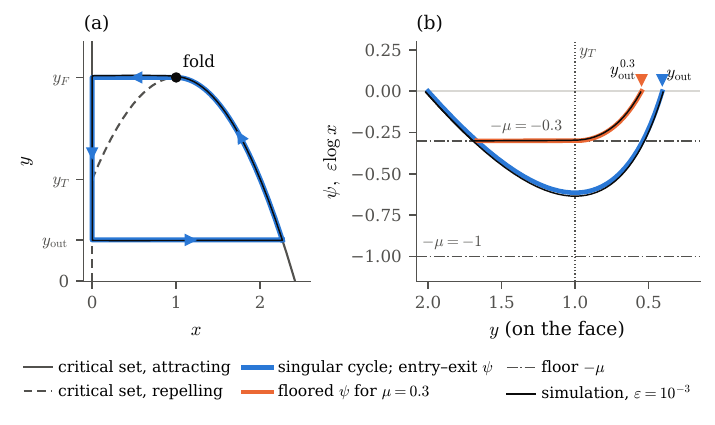}
\caption{The fold--face cycle of network $\Net1$. (a)~Phase plane. Grey: the critical
set, that is, the branch $y=1+2x-x^2$ and the face $x=0$, solid where attracting and
dashed where repelling; dot: the fold $(x_F,y_F)=(1,2)$. Thick blue: the singular cycle,
with arrows in the direction of motion; it drops from the fold to the face, follows the
face down to $y_{\rm out}=0.40638$ and rises to the branch at $x_+=2.2624$. Thin black:
the cycle for $\eps=10^{-3}$. (b)~The passage along the face on the logarithmic scale
$\psi=\eps\log x$, with $y$ decreasing from $y_F$. Thick blue: the entry--exit curve
$\psi(y)=2[(y-2)-\ln(y/2)]$, with minimum $-D=-2(1-\ln2)$ at $y_T=1$ and zero at
$y_{\rm out}$. Thick orange: the floored curve for the breaking reaction $0\to\ce X$ of
rate $e^{-\mu/\eps}$ with $\mu=0.3$, with exit $y^{0.3}_{\rm out}=0.54737$; for $\mu=1>D$
the floor is not reached and the curve is the blue one. Dash-dotted: the floors $-\mu$.
Thin black: $\eps\log x$ for $\eps=10^{-3}$, on the cycle of (a) and, for $\mu=0.3$,
started just after the drop; they lie within $0.021$ of the predicted curves.}
\label{fig:cycle}
\end{figure}

\paragraph{Example (network $\Net1$).} Fast reactions
$\ce X\to2\ce X$ $(1)$, $2\ce X\to3\ce X$ $(2)$, $3\ce X\to2\ce X$ $(1)$,
$\ce X+\ce Y\to\ce Y$ $(1)$; slow reactions $\ce X+\ce Y\to\ce X+2\ce Y$ $(c=1)$,
$\ce Y\to0$ $(m=\tfrac12)$.
\begin{itemize}[nosep]
  \item $A=1+2x-x^2$, $B=-1$, $P=2x-2$.
  \item $H=1+2x-x^2$, with $y_T=1$ and fold $(x_F,y_F)=(1,2)$.
  \item $x^*=m/c=\tfrac12<x_F$, so (N1)--(N5) hold.
  \item The exit solves $\ln(y/2)=y-2$, giving $y_{\rm out}=0.40638$. The rise
        lands at $x_+=2.2624$.
  \item Period: $T_0=4.3310$, of which $3.1872$ on the face and $1.1437$ on the
        branch.
\end{itemize}
Table~\ref{tab:ro} compares this with simulations, and Figure~\ref{fig:cycle}(a)
shows the cycle.

\subsection{Minimal fold--face networks}\label{sec:enum}

We now classify the minimal fold--face networks within a class of planar
networks; outside it there are others (Remark~\ref{rem:outside}). The universe is
that of \citet{BeldimanHalanay2026}: species $\ce X$ (fast) and $\ce Y$ (slow),
molecularity at most three. A reaction $a\ce X+b\ce Y\to a'\ce X+b'\ce Y$
($a+b,\,a'+b'\le3$) changes $\ce X$ by $N=a'-a$ and $\ce Y$ by $M=b'-b$. We use
the species scaling of that paper, $x'=f$, $y'=\eps g$ in fast time, with
\[
  f=\sum_rN_rk_rx^{a_r}y^{b_r},\qquad g=\sum_rM_rk_rx^{a_r}y^{b_r}.
\]
The universe has $90$ reactions. The \emph{class} consists of the
$55$ reactions that:
\begin{itemize}[nosep]
  \item keep $\{\ce X\}$ a siphon: no reaction with $a=0$ and $N>0$;
  \item keep $f/x$ affine in $y$ with negative slope: $N\neq0$ implies $b=0$, or
        $b=1$ and $N<0$.
\end{itemize}
Then $f=x\varphi$ with $\varphi=A(x)+yB(x)$, where
\[
  A=\alpha_0+\alpha_1x+\alpha_2x^2,\qquad B=\beta_0+\beta_1x .
\]
Here $\alpha_{a-1}$ is the net $\ce X$-change of the reactions with $b=0$ and
$\ce X$-order $a$, weighted by their rate constants, and $\beta_{a-1}$ likewise
for $b=1$. The class forces $\beta_0,\beta_1\le0$. Call a network of the class
\emph{fold--face} if (N1)--(N5) hold for some rate constants.

\begin{lemma}[Families]\label{lem:fam}
Within the class, (N1)--(N4) hold for some rate constants if and only if the network
contains:
\begin{enumerate}[label=(F\arabic*),nosep]
  \item a reaction with $a=1$, $b=0$, $N>0$: $\ce X\to2\ce X$, $\ce X\to3\ce X$ or
        $\ce X\to2\ce X+\ce Y$;
  \item a reaction with $a=1$, $b=1$, $N<0$: $\ce X+\ce Y\to0$, $\ce Y$, $2\ce Y$
        or $3\ce Y$;
  \item $2\ce X\to3\ce X$ (the only reaction with $a=2$, $b=0$, $N>0$);
  \item a reaction with $a=3$, $b=0$, $N<0$ (nine choices);
  \item an $\ce X$-free reaction consuming $\ce Y$ ($\ce Y\to0$, $2\ce Y\to\ce Y$,
        $2\ce Y\to0$, $3\ce Y\to2\ce Y$, $3\ce Y\to\ce Y$, $3\ce Y\to0$).
\end{enumerate}
In addition, no reaction may start from the complex $0$, and the lowest and the
highest $\ce Y$-orders among the $\ce X$-free reactions must each carry a consumer.
Three reactions of the class start from $0$ ($0\to\ce Y$, $2\ce Y$, $3\ce Y$).
\end{lemma}

\begin{proof}
\emph{(N1).} $\alpha_0>0$ needs (F1) and $\beta_0<0$ needs (F2).

\emph{(N2) and (N3).} The fold polynomial is
\[
  P=A'B-AB'=(\alpha_1\beta_0-\alpha_0\beta_1)+2\alpha_2\beta_0\,x+\alpha_2\beta_1\,x^2 .
\]
We have $\alpha_0\beta_1\le0$ and $\beta_0<0$. So $P(0)<0$ forces
$\alpha_1>\alpha_0\beta_1/\beta_0\ge0$, which needs (F3).

If $\alpha_2\ge0$, all coefficients of $P$ are $\le0$, so $P$ has no positive
root. Hence $\alpha_2<0$, which needs (F4). Conversely, $\alpha_2<0$ makes the
linear coefficient positive and the quadratic one nonnegative, so $P$ has exactly
one positive root, and it is simple. It also makes $A$, and therefore $H$,
negative at infinity, which is (N3). The inequality
$\alpha_1|\beta_0|>\alpha_0|\beta_1|$ is reached by a large rate of
$2\ce X\to3\ce X$.

\emph{(N4).} The $\ce X$-free reactions do not change $\ce X$, by the siphon
condition. On the face, $g(0,y)=\sum_bg_by^b$, where each $g_b$ is a sum
of rate constants with the signs of the $M$ of the reactions from $b\ce Y$. Since
$g_0\ge0$, $g(0,y)<0$ near $y=0$ needs $g_0=0$, so no reaction may start
from $0$. Negativity for all $y>0$ needs the lowest and the highest present orders
to be negative. When they are, the intermediate terms can be dominated by
scaling.
\end{proof}

\begin{lemma}[Locking certificate]\label{lem:lock}
Group the reactions with $a\ge1$ by $(a,b)$. Four groups have a net $\ce X$-change
of known sign under (N1)--(N3): $(1,0)$ and $(2,0)$ are positive ($\alpha_0>0$, and
$\alpha_1>0$ by the proof of Lemma~\ref{lem:fam}), and $(1,1)$ and $(3,0)$ are
negative ($\beta_0<0$, $\alpha_2<0$). Suppose there are constants $(c,d)$ such
that, with $m_{a}=c+d(a-1)$, there are multipliers $m'_G$ for which:
\begin{itemize}[nosep]
  \item $M_r\le m'_GN_r$ for every reaction $r$ in every group $G$;
  \item $m'_G=m_a$ for the groups of unknown sign;
  \item $m'_{(1,0)}\le m_1$, $m'_{(2,0)}\le m_2$, $m'_{(1,1)}\ge m_1$ and
        $m'_{(3,0)}\ge m_3$.
\end{itemize}
Then $g(x_F,y_F)<0$ at the fold for all rate constants satisfying (N1)--(N4), so
(N5) fails.
\end{lemma}

\begin{proof}
At the fold, $\varphi=\varphi_x=0$, so for any $c$ and $d$
\[
  g=g-c\,x\varphi-d\,x^2\varphi_x=\sum_{a_r\ge1}k_rx^{a_r}y^{b_r}\bigl(M_r-m_{a_r}N_r\bigr)+g(0,y).
\]
Take a group $G$ whose net $\ce X$-change $\nu_G=\sum_{r\in G}k_rN_r$ is positive,
and let $m'_G\le m_a$. Its part of the sum is
\[
  x^ay^b\Bigl(\sum_{r\in G}k_r(M_r-m'_GN_r)-(m_a-m'_G)\nu_G\Bigr)\le0 .
\]
The groups of negative sign are handled the same way with $m'_G\ge m_a$. The
$\ce X$-free part $g(0,y)$ is negative by (N4).
\end{proof}

The existence of a certificate is an exact linear feasibility problem in $(c,d)$,
which we decide by Fourier--Motzkin elimination in rational arithmetic. It is the
fold--face analogue of the yield-diagram test of
\citet{BeldimanHalanay2026rankone}: the points are $(a_r,M_r/N_r)$ and the lines
are $m_a=c+d(a-1)$. Adding reactions only adds constraints, so ``no certificate''
(\emph{unlocked}) is preserved. It also implies that some reaction with $a\ge1$
produces $\ce Y$, since $c=d=0$ is a certificate otherwise.
The script implements the lemma without the shift of the group $(2,0)$, a weaker
certificate. This does not affect the count below: a stronger certificate can only
lock more sets, and it cannot lock the networks of Theorem~\ref{thm:enum}, which
have witnesses.

\begin{lemma}[Helly bound]\label{lem:helly}
The certificate system consists of half-planes in $(c,d)$. Each comes from one
reaction, except a constant test that involves two reactions of one group. (A
single reaction with $N=0$ and $M>0$ gives an empty half-plane, and unlocks by
itself.) If it is infeasible, then by Helly's theorem in $\R^2$ some three of the constraints are
already infeasible. Hence a minimal unlocked network containing all families is one
representative per family plus at most three reactions: at most eight reactions.
\end{lemma}

The script \nolinkurl{julia/models/p3_enumerate.jl} (about 12 minutes; output
\nolinkurl{julia/results/p3_enumerate.txt} and
\nolinkurl{julia/results/p3_foldface_minimal.txt}, which lists each network with
its witness rate constants and singular cycle) proceeds as follows. (In its output,
family (F5) is labelled F6, and F5 denotes the reactions with $a\ge1$ that produce
$\ce Y$.)
\begin{enumerate}[nosep]
  \item It forms the $3\cdot4\cdot1\cdot9\cdot6=648$ sets with one reaction from
        each family. Of these, $450$ are unlocked and $198$ are locked.
  \item It extends the locked sets by reactions with $a\ge1$, one at a time, up to
        size eight (Lemma~\ref{lem:helly}). It keeps the unlocked sets all of
        whose proper subsets with all families are locked.
  \item It searches for witness rate constants satisfying (N1)--(N5) for every
        minimal set, drawing up to $4000$ log-normal random rate vectors per set from a
        fixed seed. The test computes the fold, the exit point $y_{\rm out}$ by the
        entry--exit integral, the landing point $x_+$, and $g$ along the branch.
\end{enumerate}

\begin{theorem}\label{thm:enum}
Every fold--face network of the class contains a minimal one, and there are exactly
$3474$ minimal fold--face networks in the class:
\begin{itemize}[nosep]
  \item $450$ with five reactions;
  \item $3024$ with six, each obtained by adding to a locked five-reaction set a
        reaction that unlocks it;
  \item none with seven or more.
\end{itemize}
\end{theorem}

\begin{proof}
A fold--face network contains all families (Lemma~\ref{lem:fam}) and is unlocked
(Lemma~\ref{lem:lock}). So it contains a minimal unlocked set with all families,
which has at most eight reactions (Lemma~\ref{lem:helly}); the enumeration lists
all of them. Every set in the list has a witness, checked in exact arithmetic below,
so it is fold--face. Its proper
subsets either lack a family or are locked, so they are not fold--face, and the set
is minimal. As a consistency check, none of the $198$ locked five-reaction sets
yielded a witness in a random search of $300$ rate vectors per set.
\end{proof}

The exclusions in this proof are exact: \nolinkurl{p3_enumerate.jl} decides locking
by Fourier--Motzkin elimination in rational arithmetic. The witnesses are found in
floating point, with (N5) required on $[x_F,1.001\,x_+]$, so that no witness lands
next to a zero of $g$. The script \nolinkurl{julia/models/p3_witness_check.jl}
(its output is \nolinkurl{julia/results/p3_witness_check.txt}) checks them again,
independently of the enumeration and in exact rational arithmetic. So
Theorem~\ref{thm:enum} is a computer-assisted result with no floating-point step:
the enumeration certifies the exclusions, and this script the witnesses. It reads
the listed rate constants as exact decimals and checks that they are positive, that
every reaction has molecularity at most three, and that the networks are distinct:
$450$ with five reactions and $3024$ with six. It then verifies (N1)--(N4) and
membership of the class. For (N5) it needs an interval with rational ends that
contains $[x_F,x_+]$. It encloses $x_F$ by bisection on $P$, and $y_F=H(x_F)$ from
above. It then certifies a rational lower bound $y_{\rm lo}$ of $y_{\rm out}$. The
integral $I(y)=\int_{y_F}^{y}\lambda/g(0,\cdot)$ is negative on $(y_{\rm out},y_F)$
and positive below $y_{\rm out}$, so it suffices that $I(y_{\rm lo})>0$. An upper
Riemann sum shows this, on a geometric partition, with each piece bounded by
interval arithmetic. The sum runs up to the upper bound of $y_F$, which only
enlarges it, since $\lambda/g(0,\cdot)>0$ above $y_T$. Since $H$ decreases beyond
$x_F$, a rational $x_\varphi>x_F$ with $H(x_\varphi)\le y_{\rm lo}$, both checked
exactly, satisfies $x_\varphi>x_+$. A Sturm sequence for the numerator of
$g(x,H(x))$ shows that it has no root on $[x_F,x_\varphi]$. All $3474$ witnesses
pass, with $y_{\rm lo}$ within $1\%$ of the floating-point value of $y_{\rm out}$.
Two sampled quantities, in floating point and not part of the certificate, show the
margins. The smallest relative margin of $g$ along the branch, the minimum of
$\sum(\text{terms})/\sum|\text{terms}|$ on a grid, is $1.5\cdot10^{-4}$. The values
of $y_{\rm out}$ listed by the enumeration, computed by Simpson's rule, agree with a
Gauss--Legendre quadrature to within $0.25\%$.

\paragraph{The smallest examples.} A five-reaction network has one reaction from
each family: autocatalysis, predation, $2\ce X\to3\ce X$, a $3\ce X$ step and an
$\ce X$-free consumption of $\ce Y$. Locking depends only on the predation step,
the autocatalytic step and the $3\ce X$ step:
\begin{itemize}[nosep]
  \item predation $\ce X+\ce Y\to2\ce Y$ or $3\ce Y$: always unlocked;
  \item predation $\ce X+\ce Y\to\ce Y$: unlocked iff the autocatalytic step is
        $\ce X\to2\ce X+\ce Y$ or the $3\ce X$ step produces $\ce Y$;
  \item predation $\ce X+\ce Y\to0$: always locked. A sixth reaction is needed.
\end{itemize}
This accounts for all $198=(27+6)\cdot6$ locked sets. A typical five-reaction
network is the predator--prey network
\[
  \ce X\to2\ce X,\quad\ce X+\ce Y\to2\ce Y,\quad2\ce X\to3\ce X,\quad3\ce X\to2\ce X,\quad\ce Y\to0:
\]
prey growth with an Allee effect, predation converting prey into predator, and
predator death. Network $\Net1$ is a six-reaction one. Its five-reaction core
$\{\ce X\to2\ce X,\ \ce X+\ce Y\to\ce Y,\ 2\ce X\to3\ce X,\ 3\ce X\to2\ce X,\ \ce Y\to0\}$
is locked with certificate $(c,d)=(0,0)$, since no reaction with $\ce X$ produces
$\ce Y$, and the catalytic $\ce X+\ce Y\to\ce X+2\ce Y$ unlocks it.

\paragraph{Why $\ce X+\ce Y\to0$ locks.} With this predation step,
$(c,d)=(1,-1)$ is a certificate for every five-reaction set:
\begin{itemize}[nosep]
  \item $\ce X+\ce Y\to0$ has $M/N=1\ge c$;
  \item each autocatalytic step has $M/N\le1=c$;
  \item $2\ce X\to3\ce X$ has $M=0\le(c+d)N=0$;
  \item each $3\ce X$ step has $M/N\ge-1=c+2d$.
\end{itemize}
In words: at the fold the predator is destroyed by predation at the rate at which
prey is, and the fold equations leave no room for a net production of $\ce Y$
unless a further reaction supplies it.

\paragraph{Checks.}
\begin{itemize}[nosep]
  \item \emph{Exactness of the certificate.} Two earlier enumerations found no set
        that lacks a witness and is not locked. A first enumeration without the
        certificate tested the family sets only, plus one $\ce Y$-producing
        reaction. It found witnesses for $1530$ of $1656$, and all $126$ failures
        are locked. A second enumeration used only the ungrouped certificate
        (one multiplier line, no shifts). It found witnesses for $3474$ of $3546$,
        and the $72$ failures, such as
        $\{\ce X\to2\ce X,\ \ce X+\ce Y\to0,\ 2\ce X\to3\ce X,\ 3\ce X\to\ce Y,\ \ce Y\to0,\ \ce X\to0\}$,
        are exactly those locked by the grouped form. There the slow speed at the
        fold is $-\alpha_0x_F-2k_4x_F^3-k_5y_F<0$ (rate constants numbered in the order
listed), which uses $\alpha_0>0$.
  \item \emph{Relation to canards.} Every minimal fold--face network contains a
        minimal canard network of \citet{BeldimanHalanay2026}, in line with
        Remark~\ref{rem:explosion}.
  \item \emph{Simulations.} Twelve minimal networks, drawn at random with their
        witness rate constants, were integrated in slow time from the landing
        point $(x_+,y_{\rm out})$ for six singular periods, with tolerances $10^{-8}$ and
        $10^{-10}$. The singular periods
        range from $1$ to $4\cdot10^4$, since the rates span several orders of
        magnitude. The relative error of the $y$-range against the singular cycle
        is $4\cdot10^{-4}$ to $6\cdot10^{-2}$ at $\eps=10^{-4}$, and decreases in
        all twelve at $\eps=10^{-5}$, to between $10^{-4}$ (or less) and
        $2.5\cdot10^{-2}$.
  \item \emph{Molecularity.} By Lemma~\ref{lem:fam}, every fold--face network of
        the class contains $2\ce X\to3\ce X$ and a reaction with three $\ce X$
        molecules among its reactants. So no bimolecular network in the class has
        a fold--face cycle, consistent with the absence of limit cycles in planar
        mass-action systems with reactant complexes of order at most two
        \citep{Pota1983}.
\end{itemize}

\begin{remark}[Outside the class]\label{rem:outside}
The class excludes, for instance, $2\ce X+\ce Y\to3\ce X$ ($b=1$, $N>0$). The
network
\[
  \ce X\to2\ce X\ (1),\quad\ce X+\ce Y\to2\ce Y\ (1),\quad2\ce X+\ce Y\to3\ce X\ (\tfrac12),\quad
  3\ce X\to2\ce X\ (1),\quad\ce Y\to0\ (\tfrac1{10})
\]
has, in the species scaling, $\varphi=1-y+\tfrac12xy-x^2$ and
$g=y\,(x-\tfrac1{10}-\tfrac12x^2)$, so $H=(1-x^2)/(1-x/2)$. It has a fold at
$x_F=2-\sqrt3$, $y_F=1.07180$, the transcritical point $y_T=1$, the exit
$y_{\rm out}=0.93148$ (the root of $\ln(y/y_F)=y-y_F$ below $1$), and the landing
point $x_+=0.58322$. Conditions (C1)--(C4) of Definition~\ref{def:ff} hold, as part~4
of \nolinkurl{p3_relaxation.jl} checks, so Theorem~\ref{thm:ro} applies; the
simulated cycles have $y\in[0.92892,1.07472]$ at $\eps=10^{-4}$. Removing any
reaction destroys the face exit, the fold, the face flow or the slow flow on the
branch, so the network is minimal. It contains none of the $3474$ networks,
since it lacks $2\ce X\to3\ce X$. Extending the enumeration beyond the class is
open.
\end{remark}

%======================================================================
\section{Siphon breaking}\label{sec:break}
%======================================================================

\subsection{Setting}

Let $\Sigma$ be a siphon of the fast species, and treat the slow concentration
$u\in\R$ as a parameter of the fast system $\dot z=F(z,u)=M(u,z,0)z$. Write
\[
  F(z,u)=\Lambda(u)z+Q_u(z,z)+O(|z|^3),
\]
where $Q_u$ is symmetric bilinear and collects the reactions with exactly two
siphon molecules among their reactants. Assume:
\begin{enumerate}[label=(T\arabic*),nosep]
  \item There is a \emph{critical component} $C^*$ with $s_{C^*}(u^*)=0$ and
        $s_{C^*}'(u^*)\neq0$. All other components have $s_C(u^*)<0$.
  \item $c_2:=\ell\cdot Q_{u^*}(r,r)\neq0$, where $r\ge0$ and $\ell\ge0$ are the
        right and left null vectors of $\Lambda(u^*)$, with $\ell\cdot r=1$.
\end{enumerate}
By (T1), $0$ is a simple eigenvalue of the block-triangular Metzler matrix
$\Lambda(u^*)$. So $r$ is supported on $C^*$ and the components downstream of it,
$r_D=-\Lambda_{DD}^{-1}\Lambda_{DC^*}r_{C^*}\ge0$ with intermediate components
included recursively, and $\ell$ on $C^*$ and the components upstream of it. Here
upstream and downstream refer to the linear fast edges of $\Lambda(u^*)$, not to the
feeding relation of Definition~\ref{def:graph}.
Then $s(\Lambda(u))$ has the simple root $u^*$, with
$s'(u^*)=\ell\,\Lambda'(u^*)\,r$.

A \emph{breaking reaction} has no siphon species among its reactants and produces
siphon species with net vector $b\ge0$, $b\neq0$. It is added to the fast
subnetwork with rate $\delta\,c(u)$, $c>0$. So the critical set near $(0,u^*)$ is
the zero set of
\[
  F(z,u,\delta)=F_0(z,u)+\delta c(u)\,b,\qquad F_0(z,u)=\Lambda(u)z+Q_u(z,z)+O(|z|^3),
\]
in the closed orthant. For one siphon species and $\delta=0$, $u^*$ is the
transcritical point at which the interior branch crosses the face
\citep{KrupaSzmolyan2001trans}.

\paragraph{Upstream and downstream.}
\begin{itemize}[nosep]
  \item Let $S=\{i:\ell_i>0\}$. It consists of $C^*$ and the components upstream of
        it. Let $T$ be its complement.
  \item For $j\in T$, $0=(\ell\Lambda)_j=\sum_{i\in S}\ell_i\Lambda_{ij}$, and the
        terms are $\ge0$. Hence $\Lambda_{ST}=0$: no linear edge leads from $T$
        into $S$.
  \item $\Lambda_{TT}(u^*)$ is Hurwitz, since all components in $T$ have $s_C<0$.
\end{itemize}

\begin{definition}[Breaking order]\label{def:order}
Let $w^T\in\{1,2,\dots,\infty\}^T$ be the greatest vector with
\begin{itemize}[nosep]
  \item $w^T_j=1$ if $b_j>0$;
  \item $w^T_i\le y_r\cdot w^T=\sum_jy_{rj}w^T_j$ for every fast reaction $r$
        with positive net production of $X_i\in T$ whose siphon reactants are
        nonempty and all lie in $T$.
\end{itemize}
Equivalently, $w^T$ is computed by min-plus iteration from $w^T_j=1$ on the support
of $b$ and $w^T_j=\infty$ elsewhere: replace $w^T_i$ by the minimum of $w^T_i$ and
$y_r\cdot w^T$ over these reactions until nothing changes (the iteration stops,
since the values decrease in a well-ordered set).
The \emph{breaking order} is
\[
  w=\begin{cases}
     1 & \text{if } b_S\ne0, \text{ i.e. } \ell\cdot b>0,\\[2pt]
     \min\{y_r\cdot w^T:\ r\in\mathcal R_{T\to S}\} & \text{otherwise},
    \end{cases}
\]
where $\mathcal R_{T\to S}$ is the set of fast reactions whose siphon reactants all
lie in $T$ and which produce a species of $S$, and $\min\emptyset=\infty$. Since
$\Lambda_{ST}=0$, the reactions in $\mathcal R_{T\to S}$ have at least two siphon
molecules among their reactants, so $w\ge2$ when $\ell\cdot b=0$.
\end{definition}

This is the min-plus counterpart of the tropical recursion \eqref{eq:trop}, with
$\delta$ in the role of $e^{-1/\eps}$. It agrees with the minimum over $C^*$ of the
closure over all species: entering $S$ through species of $S$ can only raise the
weight.

\subsection{Reduction}

\begin{lemma}[Baseline]\label{lem:baseline}
For $(u,\delta)$ near $(u^*,0)$, the equation $F_T(0_S,z_T,u,\delta)=0$ has a unique
small solution $z^0_T(u,\delta)$. It satisfies:
\begin{itemize}[nosep]
  \item $z^0_j\equiv0$ if $w^T_j=\infty$;
  \item $z^0_j=\delta^{w^T_j}\zeta_j(u,\delta)$ otherwise, with $\zeta_j$ smooth and
        $\zeta_j(u^*,0)>0$.
\end{itemize}
Write $z^0=(0_S,z^0_T)$. Then $\ell\cdot F(z^0,u,\delta)=\delta^wk_0(u,\delta)$,
with $k_0$ smooth and $k_0(u^*,0)>0$ if $w<\infty$. If $w=\infty$, then $z^0$ is a
zero of $F$.
\end{lemma}

\begin{proof}
\emph{Existence, uniqueness, and the species with $w^T=\infty$.} Since
$D_{z_T}F_T(0,u^*,0)=\Lambda_{TT}(u^*)$ is invertible, the implicit function theorem
gives a unique small solution. Let $T_\infty=\{j\in T:w^T_j=\infty\}$. A reaction
producing $X_j\in T_\infty$ either has a siphon reactant in $S$, and so vanishes at
$z_S=0$, or has one in $T_\infty$. Also $b_{T_\infty}=0$. So $\{z_{T_\infty}=0\}$ is
invariant, and solving inside it gives a small solution with $z_{T_\infty}=0$. By
uniqueness, it is $z^0_T$.

\emph{Rescaling.} For the other species put $z_j=\delta^{w_j}\zeta_j$, writing $w_j$
for $w^T_j$, and divide equation $j$ by $\delta^{w_j}$. Every term then carries a
power $\delta^e$ with $e\ge0$ an integer:
\begin{itemize}[nosep]
  \item a production term with reactants $y$ has $e=y\cdot w-w_j\ge0$, by the
        defining inequality;
  \item the breaking term has $e=1-w_j=0$, since $b_j>0$ forces $w_j=1$;
  \item a consumption term has $y_j\ge1$, so $e=y\cdot w-w_j\ge0$, and $e\ge1$ if
        $|y|\ge2$.
\end{itemize}
So the rescaled system $G(\zeta,u,\delta)=0$ is smooth.

\emph{The system at $\delta=0$.} Group the species by level $m=w_j$. At $\delta=0$
only the terms with $e=0$ survive:
\begin{itemize}[nosep]
  \item the linear terms inside level $m$, a principal submatrix $\Lambda^{(m)}$ of
        $\Lambda_{TT}$;
  \item the breaking term, at level $1$;
  \item production terms with $|y|\ge2$ and $y\cdot w=m$, whose reactants have
        lower levels.
\end{itemize}
A linear edge $X_i\to X_j$ has $w_j\le w_i$, so no linear term goes from a lower
level to a higher one. Hence the system at $\delta=0$ is block triangular by
level. Each $\Lambda^{(m)}$ is Hurwitz and Metzler, as a principal submatrix of a
Hurwitz Metzler matrix, so $(-\Lambda^{(m)})^{-1}\ge0$. Solving level by level
gives $\zeta(u,0)\ge0$.

\emph{Positivity.} Since $w^T$ is the greatest solution, the value $w_j$ is attained,
and it is attained through a source: if on a set of species of level $m$ the
minima were attained only by linear edges inside the set, raising all of them to
$m+1$ would give a larger solution.
Either $b_j>0$, or a reaction with $|y|\ge2$ and $y\cdot w=w_j$ produces $X_j$, or
$X_j$ is reached inside level $w_j$ by linear edges from such a species. The
corresponding source is positive, by induction on the level, and
$(-\Lambda^{(m)})^{-1}_{ji}>0$ along a path $i\to j$. So $\zeta(u^*,0)>0$.

\emph{Smoothness in $\delta$.} $D_\zeta G$ at $\delta=0$ is block triangular with
the invertible blocks $\Lambda^{(m)}$. The implicit function theorem gives a smooth
$\zeta(u,\delta)$, and by uniqueness $z^0_T=(\delta^{w_j}\zeta_j)_j$.

\emph{The projection.} Since $\Lambda_{ST}=0$ and $z^0_S=0$,
\[
  F_S(z^0)=\delta cb_S+\sum_rn_{r,S}\,k_r(z^0)^{y_r},
\]
summed over the reactions with siphon reactants in $T$ that produce species of $S$,
with $k_r$ the rate constant times the monomial in $u$. They consume nothing in
$S$, so $n_{r,S}\ge0$. Every term of $\ell\cdot F_S(z^0)$ is therefore
nonnegative. The smallest power of $\delta$ that occurs is $\delta^w$, with a
positive coefficient. If $w=\infty$ there is no term, and then $F(z^0)=0$.
\end{proof}

\begin{lemma}[Lyapunov--Schmidt reduction]\label{lem:LS}
Write $z=z^0(u,\delta)+\eta r+v$ with $\ell\cdot v=0$. Near $(0,u^*,0)$ there is a
smooth $V(\eta,u,\delta)$ with $\ell\cdot V=0$ and $V(0,u^*,0)=0$ such that
\[
  F=0\iff v=V\ \text{and}\ \mathfrak g(\eta,u,\delta):=\ell\cdot F(z^0+\eta r+V,u,\delta)=0 .
\]
Moreover:
\begin{enumerate}[label=(\roman*),nosep]
  \item $\hat s(u):=\mathfrak g_\eta(0,u,0)$ satisfies $\hat s(u^*)=0$ and
        $\hat s'(u^*)=s'(u^*)$. Also $\hat s=s(\Lambda(u))\,(1+O(u-u^*))$, and
        $c_2(u):=\tfrac12\mathfrak g_{\eta\eta}(0,u,0)$ has $c_2(u^*)=c_2$.
  \item $\kappa(u,\delta):=\mathfrak g(0,u,\delta)=\delta^wk(u,\delta)$, with $k$ smooth and
        $k(u^*,0)=k_0(u^*,0)$. If $w=\infty$, then $\mathfrak g(0,u,\delta)\equiv0$.
  \item $\varsigma(u,\delta):=\mathfrak g_\eta(0,u,\delta)-\hat s(u)=O(\delta)$, and
        $\varsigma(u^*,\delta)=2\ell\cdot Q_{u^*}(z^0,r)+O(\delta^2+\delta^w)$.
  \item At a zero of $F$, $D_zF$ has one eigenvalue $\mu$ near $0$. It is real,
        with $\mu=\mathfrak g_\eta\,(1+o(1))$. The others have real part $\le-c<0$.
\end{enumerate}
\end{lemma}

\begin{proof}
\emph{Reduction.} $\Lambda(u^*)$ maps $\ell^\perp$ into itself and is
injective there, because $\ell\cdot r=1$. So the equation $(I-r\ell)F=0$ can be
solved for $v\in\ell^\perp$ by the implicit function theorem. Differentiating
$(I-r\ell)F(z^0+\eta r+V)\equiv0$ in $\eta$ gives
\begin{equation}\label{eq:LSid}
  D_zF\,(r+V_\eta)=\mathfrak g_\eta\,r .
\end{equation}

\emph{(i).} At $\delta=0$, $z^0=0$ and $V(0,u,0)=0$. By~\eqref{eq:LSid},
$\hat s'(u^*)=\ell\Lambda'(u^*)r=s'(u^*)$, because $V_\eta(0,u^*,0)=0$. Also
$\mathfrak g_{\eta\eta}(0,u^*,0)=\ell\cdot2Q(r,r)$. By~(iv) at $z=0$, $\delta=0$, where the
critical eigenvalue is $s(\Lambda(u))$, $\hat s$ and $s$ have the same zero set
near $u^*$. Both are smooth, with a simple zero at $u^*$ and the same derivative
there, so $\hat s/s$ is smooth and equal to $1$ at $u^*$.

\emph{(ii).} $(I-r\ell)F(z^0)=O(\delta^w)$ by Lemma~\ref{lem:baseline}. So
$V(0,u,\delta)=O(\delta^w)$, and
\[
  \mathfrak g(0,u,\delta)=\ell\cdot F(z^0)+O(|u-u^*|\,\delta^w+\delta^{1+w})=O(\delta^w),
\]
uniformly in $u$. Hence the first $w-1$ derivatives of $\mathfrak g(0,u,\cdot)$ in $\delta$
vanish at $\delta=0$, and Hadamard's lemma gives $k$. If $w=\infty$, $z^0$ is a
zero of $F$, so $V(0,u,\delta)=0$ and $\mathfrak g(0,u,\delta)=0$.

\emph{(iii).} $D_zF(z)=\Lambda(u)+2Q_u(z,\cdot)+O(|z|^2)$ and $\ell\Lambda(u^*)=0$.
Also $V_\eta(0,u^*,\delta)=O(\delta)$ and $V(0,u^*,\delta)=O(\delta^w)$. Insert
these into~\eqref{eq:LSid}.

\emph{(iv).} Take the basis $\{r+V_\eta\}\cup\{\text{basis of }\ell^\perp\}$ in the
domain and the splitting $\R^k=\R r\oplus\ell^\perp$ in the target.
By~\eqref{eq:LSid} the matrix of $D_zF$ is block upper triangular, with diagonal
blocks $\mathfrak g_\eta$ and $D_W=(I-r\ell)D_zF|_{\ell^\perp}$. So
$\det D_zF=\beta\,\mathfrak g_\eta\det D_W$, where $\beta$ is close to $1$, because
$V_\eta$ is small. The $k-1$ eigenvalues of $D_W$, and the $k-1$ eigenvalues of
$D_zF$ away from $0$, are both close to the nonzero eigenvalues of $\Lambda(u^*)$:
\begin{itemize}[nosep]
  \item those of $C^*$ other than its Perron root, whose real parts are smaller,
        since $\Lambda_{C^*C^*}$ is irreducible;
  \item those of the other components, with $s_C<0$.
\end{itemize}
Dividing the two determinant formulas gives $\mu=\mathfrak g_\eta(1+o(1))$.
\end{proof}

\begin{lemma}[Positivity]\label{lem:pos}
Every zero of $F$ in a small neighbourhood of $(0,u^*,0)$, with $\delta\ge0$ and
$z_{C^*}\ge0$, has $z\ge0$.
\end{lemma}

\begin{proof}
Let $R$ be the complement of $C^*$. Split
\[
  F_R=\Lambda_{RR}x+\Lambda_{RC^*}z_{C^*}+\delta cb_R+P_R-D\,x,\qquad x=z_R,
\]
where $P_R$ collects the production terms with at least two siphon reactants, and
$D$ is diagonal and collects the coefficients $k z^{y-e_i}$ of the nonlinear
consumption terms. For $z\ge0$, $P_R\ge0$ and $D\ge0$. The matrix $\Lambda_{RR}-D$
is Hurwitz and Metzler, so $(D-\Lambda_{RR})^{-1}\ge0$.

Consider the map
\[
  \Phi(x)=(D-\Lambda_{RR})^{-1}\bigl(\Lambda_{RC^*}z_{C^*}+\delta cb_R+P_R\bigr).
\]
For small $z_{C^*}$ and $\delta$ it maps a small box $\{0\le x\le\varrho\}$ into
itself. By Brouwer's theorem it has a fixed point there. That fixed point solves
$F_R=0$, and the small solution of $F_R=0$ is unique by the implicit function
theorem. So it is $z_R$.
\end{proof}

\subsection{The theorem}

\begin{theorem}[Siphon breaking]\label{thm:break}
Assume (T1)--(T2), and let $\hat s$, $c_2(u)$, $\kappa$ and $\varsigma$ be as in
Lemma~\ref{lem:LS}.
\begin{enumerate}[label=(\alph*),nosep]
  \item ($\delta=0$.) Near $(0,u^*)$ the critical set in the orthant consists of the
        face $z=0$ and a branch $z=\eta r+O(\eta^2)$,
        $\eta=-\hat s(u)/c_2+O(\hat s^2)\ge0$. The branch crosses the face
        transcritically at $u^*$. On it the critical eigenvalue is
        $-\hat s(u)(1+o(1))$, and the other eigenvalues have negative real part.
        So:
        \begin{itemize}[nosep]
          \item if $c_2>0$, the branch lies where $s<0$ and is repelling, with one
                positive eigenvalue;
          \item if $c_2<0$, it lies where $s>0$ and is attracting.
        \end{itemize}
  \item ($w=\infty$.) $S$ together with $\{j\in T:w^T_j=\infty\}$ is a siphon of the
        fast network with the breaking reaction. The baseline $z^0$ lies on its
        face. Then $\mathfrak g=\eta\,(\hat s+\varsigma+c_2\eta+\dots)$, so the transcritical
        point persists, at the zero $u^*_\delta=u^*+O(\delta)$ of $\hat s+\varsigma$.
  \item ($w<\infty$, $c_2>0$.) For small $\delta>0$ the critical set has exactly one
        fold in the orthant within a fixed neighbourhood of $(0,u^*)$. It is
        nondegenerate, and with $\theta=\min(1,w/2)$:
        \[
          z_F=z^0+\eta_F\,r+O(\delta^{w/2+\theta}),\qquad
          \eta_F=\sqrt{\kappa/c_2}\,\bigl(1+O(\delta^\theta)\bigr),
        \]
        \[
          \hat s(u_F)+\varsigma(u_F,\delta)=-2\sqrt{c_2\,\kappa}\,\bigl(1+O(\delta^\theta)\bigr),
          \qquad \kappa=\kappa(u_F,\delta)\asymp\delta^w .
        \]
        In particular
        $u_F=u^*_\delta-2\sqrt{c_2k(u^*,0)}\,\delta^{w/2}/s'(u^*)\,(1+o(1))$, where
        $u^*_\delta$ is the zero of $\hat s+\varsigma$. The shift
        $u^*_\delta-u^*=O(\delta)$ is negligible for $w=1$, and at most of the same
        order for $w=2$. For $w\ge3$ it dominates if
        $\ell\cdot Q_{u^*}(z^0_1,r)\ne0$, where $z^0=\delta z^0_1+O(\delta^2)$.
        Otherwise it is $O(\delta^2)$: negligible for $w=3$, of the same order for
        $w=4$, and dominant for $w\ge5$ if its $\delta^2$ coefficient is
        nonzero.
        The other components carry the baseline $z^0$, of order $\delta$ where
        $w^T=1$, so the distance of the fold from the face is of order
        $\delta^{\min(1,w/2)}$; its critical component $z_{C^*}$ is of order
        $\delta^{w/2}$.
  \item ($w<\infty$, $c_2<0$.) Some fixed neighbourhood of $(0,u^*)$ contains no fold,
        for all small $\delta>0$.
\end{enumerate}
\end{theorem}

\begin{proof}
By Lemma~\ref{lem:LS}, the critical set near $(0,u^*)$ is $\{\mathfrak g=0\}$. A fold is a
point with $\mathfrak g=\mathfrak g_\eta=0$, $\mathfrak g_{\eta\eta}\ne0$ and $\mathfrak g_u\ne0$: there $u$ has a
nondegenerate extremum along the smooth curve $\{\mathfrak g=0\}$, and the critical
eigenvalue changes sign, by Lemma~\ref{lem:LS}(iv).

\emph{(a).} At $\delta=0$, $\mathfrak g=\eta\,(\hat s(u)+c_2(u)\eta+O(\eta^2))$. So the
nontrivial zero is $\eta=-\hat s/c_2+O(\hat s^2)$. On it,
$z_{C^*}=\eta r_{C^*}+O(\eta^2+\eta|u-u^*|)$, which is $\ge0$ if and only if
$\eta\ge0$. Lemma~\ref{lem:pos} then gives $z\ge0$. On the branch,
$\mathfrak g_\eta=\hat s+2c_2\eta+O(\eta^2)=-\hat s(1+o(1))$. The rest is
Lemma~\ref{lem:LS}(iv).

\emph{(b).} Let $W$ be the set in question. A reaction producing a species of $W$
with all siphon reactants outside $W$ would, by the definition of $w^T$ and $w$,
give that species a finite weight. The breaking reaction produces no species of
$W$. So $W$ is a siphon. By Lemma~\ref{lem:LS}(ii), $\mathfrak g(0,u,\delta)\equiv0$, and the
rest is as in~(a).

\emph{(c) and (d).} Use $\vartheta=\hat s(u)+\varsigma(u,\delta)$ instead of $u$; this is
allowed, since $\partial_u\vartheta=s'(u^*)+o(1)\ne0$. Then
\[
  \mathfrak g=\kappa+\vartheta\eta+\hat c\,\eta^2+O(\eta^3),\qquad \hat c=c_2+O(|\vartheta|+\delta),
\]
and $\kappa=\delta^wk$, with $\partial_\vartheta\kappa=O(\delta^w)$. The fold equations
are
\[
  \mathfrak g_\eta=\vartheta+2\hat c\eta+O(\eta^2)=0,\qquad \mathfrak g=0 .
\]
Eliminating $\vartheta$ gives $\kappa-\hat c\,\eta^2+O(\eta^3)=0$.
\begin{itemize}[nosep]
  \item If $c_2<0$, the left side is at least $\kappa>0$ for small $\eta$. So there
        is no fold, which proves~(d).
  \item If $c_2>0$, every fold in a small fixed neighbourhood has
        $\eta^2=\kappa/c_2\,(1+o(1))$. So $|\eta|$ and $|\vartheta|$ are of order
        $\delta^{w/2}$.
\end{itemize}
In the second case, rescale $\eta=\delta^{w/2}X$, $\vartheta=\delta^{w/2}Y$. The equations
become
\[
  Y+2c_2X=O(\delta^\theta),\qquad k(u^*,0)+XY+c_2X^2=O(\delta^\theta),
\]
smooth in $(X,Y,\delta^{1/2})$. At $\delta=0$ the solutions are
$X=\pm\sqrt{k(u^*,0)/c_2}$, $Y=-2c_2X$. The Jacobian determinant there is
$2c_2X\ne0$. So for small $\delta$ the implicit function theorem gives exactly two
folds, smooth in $\delta^{1/2}$, and these are all the folds in the neighbourhood.

The remainder $V=O(\delta^w+\eta(|u-u^*|+\delta)+\eta^2)$ is $o(\eta)$ there. So
the fold with $X<0$ has $z_{C^*}<0$ and lies outside the orthant. The fold with
$X>0$ has $z_{C^*}>0$, and $z\ge0$ by Lemma~\ref{lem:pos}. It is nondegenerate:
\begin{itemize}[nosep]
  \item $\mathfrak g_{\eta\eta}=2c_2+o(1)$;
  \item $\mathfrak g_u=(\eta+\partial_\vartheta\kappa+O(\eta^2))\,\partial_u\vartheta\ne0$, since
        $\partial_\vartheta\kappa=O(\delta^w)=o(\eta)$.
\end{itemize}
The formulas follow from $\eta_F^2=\kappa/\hat c\,(1+O(\eta))$ and $\vartheta=-2\hat c\eta+O(\eta^2)$.
\end{proof}

\begin{remark}
\begin{itemize}[nosep]
  \item For $w=1$, $\kappa=\delta c(u^*)\,\ell\cdot b\,(1+o(1))$.
  \item The function $\hat s$ is $s(\Lambda(u))$ up to a factor $1+O(u-u^*)$. That
        factor changes $u_F$ only at relative order $\delta^{w/2}$ for $w\le2$.
  \item In case RH of Table~\ref{tab:break}, $z^0=(0,\delta+O(\delta^2))$, and the
        hyperedge $2\ce X_2\to\ce X_1+2\ce X_2$ gives $Q(z,z)_1\ni x_2^2$. So
        $\varsigma(u^*,\delta)=2\ell\cdot Q(z^0,r)+O(\delta^2)=2\delta+O(\delta^2)$, the
        value used in the table.
\end{itemize}
\end{remark}

\begin{remark}[Network reading]
\begin{itemize}[nosep]
  \item $c_2=\sum_\varrho k_\varrho\,u^{y_{\varrho U}}\,\bigl(\ell\cdot\Delta_\varrho^\Sigma\bigr)\prod_jr_j^{y_{\varrho j}}$,
        summed over the fast reactions $\varrho$ with two siphon molecules among their
        reactants, where $\Delta^\Sigma_\varrho$ is the net change of the siphon
        species. It is the \emph{Perron-weighted net second-order self-interaction}
        of the siphon. Cooperative production such as $2\ce X_i\to3\ce X_i$
        counts positively, pairwise loss negatively. A hyperedge from a downstream
        species into $C^*$ counts positively, through $r_D>0$.
  \item The fold is created by the breaking reaction, but its existence depends only
        on the siphon, through $c_2$. The breaking reaction sets its size, through
        $w$ and $\kappa$.
  \item A breaking reaction that produces only species downstream of $C^*$, and
        that does not come back by hyperedges, breaks $\Sigma$ but not the smaller
        siphon around $C^*$ (case R$'$ below).
\end{itemize}
\end{remark}

\subsection{Link to the fold--face networks}

In the planar class of Section~\ref{sec:planar} ($\Sigma=\{\ce X\}$,
$\varphi=A(x)+yB(x)$, $B<0$), $\Lambda=\varphi(0,y)$, $\ell=r=1$, and the
transcritical point is at $y_T=-A(0)/B(0)$. There
$c_2=\partial_x\varphi(0,y_T)=\alpha_1+\beta_1y_T$, where
$A=\alpha_0+\alpha_1x+\dots$ and $B=\beta_0+\beta_1x+\dots$. Hence
\[
  P(0)=\alpha_1\beta_0-\alpha_0\beta_1=\beta_0\,(\alpha_1+\beta_1y_T)=\beta_0\,c_2,
\]
and since $\beta_0<0$, the condition $P(0)<0$ of (N2) is exactly $c_2>0$. Since
$H'(0)=-P(0)/B(0)^2$, it says that the interior branch leaves the transcritical
point increasing. It is only the first half of (N2), which also asks for a fold
$x_F$ further on. So:

\begin{corollary}\label{cor:planarbreak}
In the planar affine class, let $A(0)>0$, $B(0)<0$ and $H'(0)\ne0$ (that is,
$c_2\ne0$). The interior branch leaves the transcritical point $(0,y_T)$ increasing,
$H'(0)>0$, if and only if breaking the siphon $\{\ce X\}$, by an inflow or an
uncatalysed production of $\ce X$, creates a fold near $(0,y_T)$.
\end{corollary}

For example, $\ce X\to2\ce X$, $2\ce X\to3\ce X$, $\ce X+\ce Y\to\ce Y$ has
$H=1+x$: the branch rises, breaking creates a fold near $(0,1)$, but $H$ has no fold
of its own. The assumption $H'(0)\ne0$ is needed: for $\varphi=1-y+x^2$
($\ce X\to2\ce X$, $3\ce X\to4\ce X$, $\ce X+\ce Y\to\ce Y$, molecularity four),
$H=1+x^2$ has $H'(0)=0$, and an inflow of rate $\delta$ gives the critical curve
$y=1+x^2+\delta/x$, with a fold at $x=(\delta/2)^{1/3}$.

\subsection{Numerical check}\label{sec:breaknum}

The script \nolinkurl{julia/models/p3_breaking.jl} (output
\nolinkurl{julia/results/p3_breaking.txt}, about 10 seconds) computes the fold. It
solves $F=0$, $\det D_zF=0$ by Newton's method in 256-bit arithmetic, with $D_zF$
exact inside the determinant and a central finite-difference Newton matrix, and
stops when the relative step is below $10^{-45}$. The siphon is
$\{\ce X_1,\ce X_2\}$, with parameter $u$; all networks contain
\[
  \ce U+\ce X_1\to\ce U+2\ce X_1,\quad\ce X_1\to0\ (2),\quad\ce X_1\to\ce X_1+\ce X_2,\quad
  \ce X_2\to0,\quad3\ce X_1\to2\ce X_1,\quad2\ce X_2\to\ce X_2 .
\]
\begin{itemize}[nosep]
  \item \textbf{Irreducible} (I): add $\ce X_2\to\ce X_1+\ce X_2$ $(\tfrac12)$.
        Then $\Lambda=\bigl(\begin{smallmatrix}u-2&1/2\\1&-1\end{smallmatrix}\bigr)$,
        $u^*=\tfrac32$, $r=(1,1)$, $\ell=(\tfrac23,\tfrac13)$, $s'(u^*)=\tfrac23$,
        and $c_2=(2a-1)/3$ with $a$ the rate of $2\ce X_1\to3\ce X_1$.
  \item \textbf{Reducible} (R): $\ce X_1$ feeds $\ce X_2$ only. Then
        $\Lambda=\bigl(\begin{smallmatrix}u-2&0\\1&-1\end{smallmatrix}\bigr)$,
        $u^*=2$, $r=(1,1)$, $\ell=(1,0)$, $c_2=a$, plus $1$ for the hyperedge
        $2\ce X_2\to\ce X_1+2\ce X_2$ in RH.
\end{itemize}

\begin{table}[ht]
\centering
\footnotesize
\setlength{\tabcolsep}{4pt}
\begin{tabular}{@{}lllllll@{}}
\toprule
case & breaking & $c_2$ & $w$ & $\delta$ & $x_1/\eta_F$ & $(u^*-u_F)/(\text{predicted})$\\
\midrule
I$+$ ($a=2$) & $0\to\ce X_1$ & $1$ & $1$ & $10^{-4,-6,-8}$ & $1.00016,\ 1.00000,\ 1.00000$ & $0.99995,\ 1.00000,\ 1.00000$\\
I$+'$ & $0\to\ce X_2$ & $1$ & $1$ & same & $0.99994,\ 1.00000,\ 1.00000$ & $0.99425,\ 0.99942,\ 0.99994$\\
I$-$ ($a=\tfrac14$) & $0\to\ce X_1$ & $-\tfrac16$ & $1$ & same & \multicolumn{2}{l}{no fold (predicted: none)}\\
R ($a=2$) & $0\to\ce X_1$ & $2$ & $1$ & same & $1.00357,\ 1.00035,\ 1.00004$ & $0.99823,\ 0.99982,\ 0.99998$\\
R$'$ & $0\to\ce X_2$ & $2$ & $\infty$ & same & \multicolumn{2}{l}{no fold; $\{\ce X_1\}$ stays a siphon}\\
RH & $0\to\ce X_2$ & $3$ & $2$ & $10^{-3,-4,-5}$ & $1.00057,\ 1.00006,\ 1.00001$ & $0.99746,\ 0.99975,\ 0.99997$\\
\bottomrule
\end{tabular}
\caption{Siphon breaking with two siphon species: the fold $(x_1,u_F)$ computed for
each $\delta$ (one value per $\delta$ in each cell), against Theorem~\ref{thm:break}.
$\eta_F=\sqrt{\kappa/c_2}$ with
$\kappa=(\ell\cdot b)\delta$ for $w=1$ and $\kappa=\delta^2$ for RH. In RH the
predicted $u_F$ includes $\varsigma=2\delta$. Without it the last column would be
$1.577=(2+2\sqrt3)/(2\sqrt3)$, since the baseline $x_2\approx\delta$ enters the
hyperedge linearly in $\eta$.}
\label{tab:break}
\end{table}

All predictions of Table~\ref{tab:break} hold. For $c_2<0$ (I$-$) and for
$w=\infty$ (R$'$), where no fold is predicted, Newton's method from nine starting points
near $(0,u^*)$ finds none with $x>0$ within $0.2$ of $u^*$. In RH the fold is at scale $\delta$, not $\sqrt\delta$, as the
breaking order $w=2$ predicts. There the second species also carries the baseline:
$x_2=\delta+\eta_F+\dots$ ($1.577\cdot10^{-4}$ against
$\eta_F=0.577\cdot10^{-4}$ at $\delta=10^{-4}$).

%======================================================================
\section{Exponentially weak breaking: the floor}\label{sec:floor}
%======================================================================

After breaking, the critical manifold of a fold--face network has two folds: the
old one at $y_F$, and a new one at distance $O(\sqrt\delta)$ from the face, next to
$y_T$ (Theorem~\ref{thm:break}(c) and Corollary~\ref{cor:planarbreak}). For fixed
$\delta$ and $\eps\to0$ the cycle is expected to be an ordinary relaxation
oscillation on an S-shaped critical manifold \citep{KrupaSzmolyan2001relax}; we do
not check its hypotheses here. For $\delta=0$ it is the fold--face cycle.

What decides between the two is $\mu=\eps\log(1/\delta)$. On the logarithmic scale
of Section~\ref{sec:ee} the broken face is a floor: the breaking reaction is a feed
of constant logarithmic size $\eps\log\delta=-\mu$. Mutations with probability
$K^{-1/\alpha}$ act in the same way in the adaptive-dynamics limits of
Remark~\ref{rem:tropicalrelated}.

\begin{proposition}[Floor]\label{prop:floor}
Assume the setting of Theorem~\ref{thm:tropical} for the unbroken network: (S1),
(S2), acyclic feeding, and $s(0)<0$. Add a fast breaking reaction of rate $\delta c(u)$, with
$\delta=e^{-\mu/\eps}$ and $\mu>0$ fixed, that produces the siphon species $b\ge0$
and changes $u$ by at most $O(\delta)$.
Replace $F_C$ by
\[
  F^\mu_C=\max(F_C,-\mu)\quad\text{if }b_C\ne0,\qquad F^\mu_C=F_C\quad\text{otherwise},
\]
and let $\psi^\mu$ be the corresponding solution of the tropical recursion, with
exit time $\tau^\mu_{\rm ex}$. Then the conclusions of Theorem~\ref{thm:tropical}
hold with $\psi^\mu$ and $\tau^\mu_{\rm ex}$.
\end{proposition}

\begin{proof}
In the proof of Theorem~\ref{thm:tropical} the breaking reaction adds the term
$\delta c(u)\,\ell_C\cdot b_C\ge0$ to the source $R_C$. Its value satisfies
$\eps\log(\delta c\,\ell_C\cdot b_C)=-\mu+O(\eps)$ uniformly. Nothing else changes:
\begin{itemize}[nosep]
  \item the orthant is still forward invariant;
  \item the term goes into $q_C$ in the splitting \eqref{eq:metzlersplit};
  \item its effect on $u'$ is $O(e^{-\mu/\eps}/\eps)$.
\end{itemize}
Laplace's principle turns $\frac1\eps\int_0^\tau e^{\Psi_C(\sigma,\tau)/\eps}\,\delta c\,d\sigma$
into $\sup_\sigma[-\mu+\Psi_C(\sigma,\tau)]$ on the logarithmic scale, which is the
extra term in $F^\mu_C$.
\end{proof}

For one siphon species this is a reflected (Skorokhod) version of the entry--exit
relation:
\[
  \psi(\tau)=\max\Bigl(\int_0^\tau\lambda,\ \sup_{\sigma\le\tau}\Bigl[-\mu+\int_\sigma^\tau\lambda\Bigr]\Bigr),
\]
that is, $\psi'=\lambda$ except that $\psi$ stays at $-\mu$ while $\lambda<0$ there.

\begin{corollary}[Broken fold--face cycles]\label{cor:floor}
Let $\Gamma$ be the fold--face cycle of a planar network satisfying (N1)--(N5), and
add a reaction $0\to\ce X$, or any reaction producing $\ce X$ with no $\ce X$ among
its reactants, at rate $\delta c(y)$, where $c$ is smooth and positive for $y>0$ and
$\delta=e^{-\mu/\eps}$ with $\mu>0$ fixed. Let
\[
  D=-\int_{y_F}^{y_T}\frac{\lambda(y)}{g(0,y)}\,dy>0
\]
be the depth of the passage, and let $y^\mu_{\rm out}$ be given by
\begin{itemize}[nosep]
  \item $y^\mu_{\rm out}=y_{\rm out}$ if $\mu\ge D$;
  \item $y^\mu_{\rm out}\in(y_{\rm out},y_T)$ with
        $\int_{y_T}^{y^\mu_{\rm out}}\lambda(y)/g(0,y)\,dy=\mu$ if $\mu<D$.
\end{itemize}
Let $\Gamma^\mu$ be the singular cycle made of the branch from $x^\mu_+$ to the fold,
the drop, the face segment from $y_F$ down to $y^\mu_{\rm out}$, and the rise to
$x^\mu_+$, where $x^\mu_+>x_F$ and $H(x^\mu_+)=y^\mu_{\rm out}$. There is
$\eps_0(\mu)>0$ such that for $0<\eps<\eps_0(\mu)$ the network has a unique periodic
orbit in a neighbourhood of $\Gamma^\mu$. It is orbitally asymptotically stable,
with nontrivial Floquet multiplier $O(e^{-c(\mu)/\eps})$, and it tends to
$\Gamma^\mu$ in the Hausdorff distance. Its lower
turning point tends to $y^\mu_{\rm out}$, which decreases continuously from $y_T$
(as $\mu\to0$) to $y_{\rm out}$ (at $\mu=D$).
\end{corollary}

\begin{proof}
Fix $\mu>0$; constants below may depend on $\mu$. The broken system reads
\[
  \eps x'=x\varphi(x,y)+\delta c(y),\qquad y'=g(x,y)+O(\delta),
\]
where the $O(\delta)$ term, present if the breaking reaction consumes $\ce Y$, does not
depend on $x$.

\emph{The singular cycle.} The Skorokhod formula after Proposition~\ref{prop:floor},
written in the variable $y$ with $d\tau=dy/g(0,y)$, gives the exit point: $\psi^\mu$
decreases from $0$ while $y>y_T$, reaching $-D$ at $y_T$ if there is no floor. If
$\mu\ge D$ the floor plays no role and $y^\mu_{\rm out}=y_{\rm out}$. If $\mu<D$,
$\psi^\mu$ reaches $-\mu$ before $y_T$, stays there until $y_T$, where $\lambda$
changes sign, and then grows as $-\mu+\int_{y_T}^y\lambda/g(0,\cdot)$; it returns to
$0$ at $y^\mu_{\rm out}$. In both cases $y_{\rm out}\le y^\mu_{\rm out}<y_T$, so
$\lambda(y^\mu_{\rm out})>0$ and the exit is transversal. Since $H$ decreases from
$y_F$ to $y_{\rm out}$ on $[x_F,x_+]$, there is a unique $x^\mu_+\in(x_F,x_+]$ with
$H(x^\mu_+)=y^\mu_{\rm out}$. The argument for the rise in the proof of
Proposition~\ref{prop:planar}, with $y^\mu_{\rm out}$ in place of $y_{\rm out}$,
shows that the fibre at $y^\mu_{\rm out}$ rises to $x^\mu_+$, and by (N5) the slow
flow carries the branch from $x^\mu_+$ to the fold. The branch near the fold, the
fold and the drop are those of $\Gamma$. So $\Gamma^\mu$ has the four pieces of a
fold--face cycle, with $u_{\rm out}$ replaced by $y^\mu_{\rm out}$.

\emph{Away from the face.} The function $\delta=e^{-\mu/\eps}$, extended by $0$ at
$\eps=0$, is $C^\infty$ in $\eps\ge0$ and vanishes there to all orders. So the broken
system is a $C^\infty$ perturbation of the unbroken one in $(x,y,\eps)$, with the
same critical manifold at $\eps=0$. We take the sections $\Delta$ and $\Delta_3$ of
the proof of Theorem~\ref{thm:ro}, and the ball $B\subset\Delta$, on the branch
between $y^\mu_{\rm out}$ and $y_F$; they depend on $\mu$. For $\mu\ge D$, the
branch extends beyond $x_+$ as required in (C1), since $g(x,H(x))>0$ slightly beyond
$x_+$ by continuity. The maps $\Pi_1$ (branch and fold), $\Pi_2$ (drop) and $\Pi_4$
(rise and branch) are built from Fenichel's theorem, the fold theorem of
\citet{KrupaSzmolyan2001fold} and regular perturbation over bounded fast times, away
from the face or at the fixed distance $\rho$ from it. They keep the properties used
there, with the rise now starting near $y^\mu_{\rm out}$. The breaking also creates a
second fold of the critical set, at distance $O(e^{-\mu/(2\eps)})$ from the face near
$y_T$ (Theorem~\ref{thm:break}(c), with $w=1$). It needs no separate passage: the
bounds below use the full $\varphi(x,y)$ along the trajectory, so they hold through
the region of that fold.

\emph{Face passage.} Let $\mathcal K_E$ be the entry set of the proof of
Theorem~\ref{thm:ro}, a compact piece of $\{x=\rho\}$ near $(\rho,y_F)$, and let
$\tau_\eps$ be the first time $\tau\ge\tau_0$ with $x=\rho$, as in
Lemma~\ref{lem:face}. Fix $T>\tau^\mu_{\rm ex}$; on $[0,\min(\tau_\eps,T)]$ the
value of $y$ stays in a compact subset of $(0,\infty)$, and we work there until the
exit window shows $\tau_\eps<T$. The constants are chosen in this order: $\eta$, then
$m=m(\eta,\mu)$ and $s_0$ below, then $\rho\le\min(m,s_0,\rho_1)/(4K)$, and $\eps$
last. Along a trajectory put
$\Phi(\sigma,\tau)=\int_\sigma^\tau\varphi(x,y)$. By variation of constants,
\[
  x(\tau)=x(0)e^{\Phi(0,\tau)/\eps}
  +\frac1\eps\int_0^\tau e^{\Phi(\sigma,\tau)/\eps}\,\delta c(y(\sigma))\,d\sigma ,
\]
and both terms are positive. With $0<c_0\le c\le c_1$ and $|\varphi|\le K$ near the
face, this gives, for $\eps\le\tau\le\min(\tau_\eps,T)$ and uniformly on
$\mathcal K_E$,
\[
  M(\tau)-K\eps\le\eps\log x(\tau)\le M(\tau)+K\eps\log\tfrac1\eps,\qquad
  M(\tau)=\max\Bigl(\Phi(0,\tau),\ \sup_{0\le\sigma\le\tau}\bigl[-\mu+\Phi(\sigma,\tau)\bigr]\Bigr).
\]
For the lower bound, keep the first term, or the part of the integral over an
interval of length $\eps$ in $[0,\tau]$ that contains a maximiser $\sigma_0$; on it
$\Phi(\sigma,\tau)\ge\Phi(\sigma_0,\tau)-K\eps$. Write
$\varphi=\lambda(y_0)+O(x)+O(|y-y_0|)$, where $y_0$ is the face flow from the entry
value of $y$. We follow the windows of the proof of Lemma~\ref{lem:integral}.
\begin{itemize}[nosep]
  \item \emph{Entry window $[0,\tau_1]$.} Here
        $\eps x'\le x(\lambda(y)+Kx)+\delta c_1\le-\tfrac c2x+\delta c_1$, so
        $x\le\rho e^{-c\tau/(2\eps)}+K\delta$, and $\int_0^{\tau_1}x=O(\eps)$.
  \item \emph{Middle window $[\tau_1,\tau^\mu_{\rm ex}-\eta]$.} The function
        $\psi^\mu$ is negative on $(0,\tau^\mu_{\rm ex})$, so $\psi^\mu\le-m$ on the
        window for some $m\in(0,\mu]$. Since $x\le\rho$ before
        $\tau_\eps$, $|y-y_0|\le K\rho$ and so $M\le\psi^\mu+K\rho\le\psi^\mu+m/4$;
        with the upper bound above, $x\le e^{-m/(2\eps)}$ on the window. Then
        $\int x=O(\eps)$ up to $\tau^\mu_{\rm ex}-\eta$, hence $y=y_0+O(\eps)$,
        $M=\psi^\mu+O(\eps)$, and $\eps\log x=\psi^\mu+O(\eps\log\frac1\eps)$
        there. This proves Proposition~\ref{prop:floor} directly in the present
        case. On the floor $\psi^\mu=-\mu$, so $x$ is exponentially small there and
        adds nothing to $\int x$.
  \item \emph{Exit window $[\tau^\mu_{\rm ex}-\eta,\tau_\eps]$.} Here
        $\lambda(y_0)\ge s_0>0$, since $\lambda(y^\mu_{\rm out})>0$. The source is
        nonnegative, so $\eps x'\ge(s_0-Kx-K\eps)x\ge\tfrac{s_0}2x$. As in
        Lemma~\ref{lem:integral}, $x$ reaches $\rho$ at
        $\tau_\eps=\tau^\mu_{\rm ex}+o(1)$, the integral over the window is
        $O(\eps)$, and $\{x=\rho\}$ is crossed transversally.
\end{itemize}
So Lemma~\ref{lem:integral}(a)--(c) hold with $\tau^\mu_{\rm ex}$ in place of
$\tau_{\rm ex}$, uniformly on $\mathcal K_E$, and the passage map $P_\eps$ from
$\mathcal K_E$ to $\{x=\rho\}$ satisfies Lemma~\ref{lem:face}(a), with the exit
point $y^\mu_{\rm out}(y_E)$ of the entry value $y_E$ in place of $u_{\rm out}$. This
exit point is $y_{\rm out}(y_E)$ if the depth $-\int_{y_E}^{y_T}\lambda/g(0,\cdot)$
is at most $\mu$, and $y^\mu_{\rm out}$ otherwise. It is Lipschitz in $y_E$: constant
near $y_F$ if $\mu<D$, smooth if $\mu>D$, and only Lipschitz if $\mu=D$. In the
definition of $\mathcal X$, $L$ is a Lipschitz constant of this map.

\emph{Derivative of the passage.} In $\xi=\log x$, $\eps\xi'=\varphi(x,y)+\Theta$
with $\Theta=\delta c(y)/x$. The term $\Theta$ is bounded: by the argument for feeds
in Lemma~\ref{lem:apriori}(b),
$\eps(\log\Theta)'=\eps(c'/c)\,y'-\eps\xi'\le K-\Theta$, and $\Theta(0)\le K\delta$ on
$\mathcal K_E$, so $\Theta\le K$. Hence
\[
  \partial_\xi(\varphi+\Theta)=x\varphi_x-\Theta\le Kx,\qquad
  \partial_y(\varphi+\Theta)=\varphi_y+\Theta\,c'/c=O(1),\qquad
  \partial_\xi y'=xg_x=O(x).
\]
The Gronwall argument in the proof of Lemma~\ref{lem:face}(b) uses only the
one-sided bound on the logarithmic norm, and $-\Theta\le0$ only helps there, although
$\Theta$ is not $O(x)$ on the floor. With $\int_0^{\tau_\eps}x\le K\eps$,
the Gronwall argument there gives $\|DP_\eps\|\le K/\eps$.

\emph{Conclusion.} With these properties of $\Pi_1$, $\Pi_2$, $\Pi_3=P_\eps$ and
$\Pi_4$, the return map on the ball $B$ satisfies $\|D\Pi\|\le e^{-c'/\eps}$ and
maps $B$ into itself, as in the proof of Theorem~\ref{thm:ro}. Its unique fixed point
gives the unique periodic orbit in a thin neighbourhood of $\Gamma^\mu$, with
nontrivial Floquet multiplier $O(e^{-c'/\eps})$. The orbit tends to $\Gamma^\mu$
piece by piece, as there. Its lower turning point is the value of $y$ at the exit
from the face, up to $o(1)$, and so tends to $y_0(\tau^\mu_{\rm ex})=y^\mu_{\rm out}$.
Finally, $y^\mu_{\rm out}$ is continuous and decreasing in $\mu$ on $(0,D]$, tends to
$y_T$ as $\mu\to0$, and equals $y_{\rm out}$ at $\mu=D$, since
$\int_{y_T}^{y}\lambda/g(0,\cdot)$ increases from $0$ to $D$ as $y$ decreases from
$y_T$ to $y_{\rm out}$.
\end{proof}

The constants $\eps_0(\mu)$ and $c(\mu)$ tend to $0$ as $\mu\to0$, since
$m\le\mu$ and $\lambda(y^\mu_{\rm out})\to0$. So the corollary says nothing
uniform about $\mu=\mu(\eps)\to0$, the regime of fixed $\delta$ discussed next.

For fixed $\delta>0$, i.e.\ $\mu=0$ in the limit, we expect the cycle instead to
jump at the new fold near $y_T$ (Theorem~\ref{thm:break}(c)), as for an S-shaped critical manifold.
The new fold does not set the jump when $\mu>0$: past it, $x$ still has to grow from
$e^{-\mu/\eps}$ to order one, which takes the slow time in which
$\int\lambda\,d\tau=\mu$.

\paragraph{Numerical check.} Part~2 of \nolinkurl{p3_breaking.jl} checks
Corollary~\ref{cor:floor} on network $\Net1$ with $0\to\ce X$ at rate
$\delta=e^{-\mu/\eps}$. The depth is $D=2(1-\ln2)=0.61371$. Table~\ref{tab:floor} compares the
predicted and simulated lower turning points.
Figure~\ref{fig:cycle}(b) compares the floored curve for $\mu=0.3$ with $\eps\log x$
along the face for $\eps=10^{-3}$ (script \nolinkurl{p3_figures.jl}).

\begin{table}[ht]
\centering
\small
\begin{tabular}{@{}lllll@{}}
\toprule
$\mu$ & predicted & $\eps=10^{-2}$ & $\eps=10^{-3}$ & $\eps=10^{-4}$\\
\midrule
0.05 & 0.79274 & 1.11469 & 0.80791 & 0.79441\\
0.15 & 0.66100 & 0.73249 & 0.66817 & 0.66182\\
0.30 & 0.54737 & 0.58648 & 0.55169 & 0.54787\\
0.45 & 0.47002 & 0.49759 & 0.47315 & 0.47038\\
0.60 & 0.41111 & 0.43237 & 0.41356 & 0.41140\\
1.00 & 0.40638 & 0.38243 & 0.40004 & 0.40485\\
\bottomrule
\end{tabular}
\caption{Lower turning point of the cycle after breaking the siphon of $\Net1$ at
rate $e^{-\mu/\eps}$. For $\mu\ge D=0.614$ it is the entry--exit point $0.40638$,
and as $\mu\to0$ it tends to $y_T=1$.}
\label{tab:floor}
\end{table}

So for planar networks the cycles form a family parametrised by
$\mu=\eps\log(1/\delta)$, with the boundary passage at one end; at the other end we
expect the S-shaped relaxation oscillation of fixed $\delta$. A siphon broken by a
reaction whose rate is exponentially small in $1/\eps$ still behaves like a siphon
for part of the passage.

%======================================================================
\section{Discussion}\label{sec:disc}
%======================================================================

\paragraph{Relation to the literature.}
\begin{itemize}[nosep]
  \item \citet{DeMaesschalckSchecter2016} and most of the classical entry--exit
        literature treat one transverse variable.
  \item \citet{HsuRuan2021} allow several transverse variables, but each vanishes
        on its own invariant face, so the linearisation is diagonal (their
        Remark~2.1). In our terms, every component is a single species and there
        are no feeds.
  \item \citet{Borsotti2026} treats block-diagonal fast parts with nilpotent
        corrections. The exit is where the maximum of the block integrals vanishes,
        the case of \eqref{eq:trop} without feeds, with $\eps|\log\eps|$
        corrections from Jordan blocks.
  \item \citet{KaklamanosKuehnPopovicSensi2025} study two fast variables whose real
        eigenvalues intersect, in polar coordinates. Their case with an invariant
        eigendirection is a two-node instance of the tropical recursion.
  \item \citet{BacaerLobrySari2020} prove Theorem~\ref{thm:perron} for $k=2$ when the
        second eigenvalue stays negative, and
        \citet{KaklamanosPuglieseSensiSottile2024} use the principal eigenvalue of
        a $2\times2$ Metzler block in an SIRS model with secondary infections.
  \item In stochastic adaptive dynamics with mutation probability $K^{-1/\alpha}$,
        \citet{CoquilleKrautSmadi2021} and \citet{EsserKraut2025} obtain a max-plus
        recursion for the exponents of single traits fed by mutation, with rates
        that are constant between invasions or averaged over fast phases
        (Remark~\ref{rem:tropicalrelated}).
  \item In the tropical reduction of reaction networks
        \citep{RadulescuVakulenkoGrigoriev2015,Kruff2021} the orders of magnitude
        of the concentrations solve max-plus equations, the tropical
        equilibrations, and determine the time scales of the reduced models.
        Here the orders $\psi_C$ are not equilibrium values: they evolve in slow
        time along an invariant face that may be repelling.
\end{itemize}
Delayed loss of stability in several dimensions was proved by
\citet{Neishtadt1987} for a slow passage through a Hopf bifurcation in analytic
systems. For $k\ge3$, the relation with the leading eigenvalue is stated in the
literature when only that eigenvalue changes sign
\citep[see the references in][]{KaklamanosPuglieseSensiSottile2024}; we found no
statement for $k\ge3$ of Theorem~\ref{thm:perron}, which assumes nothing about the
other eigenvalues. Max-plus
recursions with feeding between types occur in the adaptive-dynamics limits above;
we found none with irreducible blocks, rates that vary along a slow flow, or
hyperedges, nor one for a deterministic passage near an invariant face.

\begin{remark}[Landing and fast non-siphon species]\label{rem:landing}
Assumption (S2) excludes fast reactions that change other species. In the setting
of \citet{Hsu2019} (the critical manifold is the face itself, and the predator,
outside the siphon, moves fast), the fast flow carries $u$ along a heteroclinic
orbit while $z$ is still of order one. The relations above should then hold from
the \emph{landing point}, the $\omega$-limit of that orbit on the face, as in
\citet{HsuRuan2021}. Fast reactions that change species outside $\Sigma$, such as
$\ce A+\ce B+\ce Y\to3\ce X$ in Olsen's model, vanish on the face, but act on $u$
during the drop and the rise. Neither extension is written out here.
\end{remark}

\begin{remark}[Entry values and memory]
Entry values $a_C<0$ occur when a component is exponentially small on arrival, for
example when it never recovered during the previous excursion away from the face.
Then \eqref{eq:trop} carries information from one passage to the next, and a
periodic orbit is a fixed point of a tropical map on the entry values.
\end{remark}

\paragraph{Open problems.}
\begin{itemize}[nosep]
  \item \emph{Cyclic feeding} (Conjecture~\ref{conj:cyclic}). It is supported by
        networks R3, R4, $\Net2$ and Olsen's model. A proof would extend
        Theorem~\ref{thm:ro} and Proposition~\ref{prop:floor} to cyclic feeding;
        the induction behind Lemma~\ref{lem:apriori} and the weights in
        Lemma~\ref{lem:face} both fail there.
  \item \emph{Error order.} Whether the error in Theorem~\ref{thm:tropical} is
        $O(\eps\log(1/\eps))$ in general.
  \item \emph{Landing cycles} and fast non-siphon species
        (Remark~\ref{rem:landing}), and cycles with memory.
  \item \emph{Canard explosions ending in face cycles}
        (Remark~\ref{rem:explosion}).
  \item \emph{Breaking order $w\ge3$.} Theorem~\ref{thm:break} predicts that the
        baseline shift $\varsigma$ then dominates the position of the fold, under a
        nondegeneracy condition; this has not been tested numerically.
  \item \emph{Planar networks outside the affine class}
        (Remark~\ref{rem:outside}).
  \item \emph{Codimension two}: several critical components at the same $u^*$.
\end{itemize}

\paragraph{Data and code availability.} The Julia scripts
\nolinkurl{julia/models/p3_*.jl} produce all tables, the list of minimal
fold--face networks and its independent check, and (\nolinkurl{p3_figures.jl}) the
data of the figures, which \nolinkurl{figures/fig_p3_tropical.py} and
\nolinkurl{figures/fig_p3_cycle.py} draw. The script \nolinkurl{p3_enumerate.jl}
also loads \nolinkurl{julia/models/planar_core.jl}, which defines the reaction universe,
and the list \nolinkurl{julia/results/planar_minimal_final.txt} of minimal canard
networks of \citet{BeldimanHalanay2026}. The scripts, their outputs
(\nolinkurl{julia/results/p3_*.txt}) and the package versions
(\nolinkurl{julia/Manifest.toml}) are openly available on Zenodo
\citep{SiphonData2026}, \url{https://doi.org/10.5281/zenodo.24259972} (data CC~BY~4.0, code MIT).

\paragraph{Use of AI tools.}
The authors used the AI assistant Claude (Anthropic), through Claude Code,
throughout this work: to write and run the Julia and Python code for the
computations and figures,
to search the literature and check the references against Crossref and
zbMATH, to draft and check the text, including the proofs, and to review the
manuscript critically. The authors directed the research, checked all
statements, proofs and computations, and take full responsibility for the
content of the paper.

\bibliographystyle{plainnat}
\bibliography{refs}

\end{document}